\documentstyle[amscd,amstex,12pt]{amsart}

\newtheorem{lemma}{Lemma}
\newtheorem{cor}{Corollary}
\newtheorem{prop}{Proposition}
\newtheorem{unnumbered}{Theorem} 
\newtheorem{theorem}{Theorem}

\theoremstyle{definition}
\newtheorem{defn}{Definition}	

\newtheorem{question}{Question}

\theoremstyle{remark}
	
\newtheorem{rems}{Remarks}	

\newcommand{\A}{\cal A}

\newcommand{\Q}{\Bbb Q}
\newcommand{\Z}{{\Bbb Z}}
\newcommand{\R}{{\Bbb R}}

\newcommand{\Gal}{\operatorname{Gal}}

\newcommand{\Hom}{\operatorname{Hom}}
\newcommand{\Res}{\operatorname{Res}}

\newcommand{\divv}{\operatorname{div}}
\newcommand{\Tr}{\operatorname{Tr}}
\newcommand{\sep}{^{\text{sep}}}

\newcommand{\Lbar}{\overline{L}}
\newcommand{\Fp}{{{\Bbb F}_p}}
\newcommand{\Fq}{{{\Bbb F}_q}}
\newcommand{\Fqq}{{{\Bbb F}_{q'}}}
\newcommand{\Fqhat}{{\hat{\Bbb F}_q}}
\newcommand{\Fpw}{{{\Bbb F}^\omega_p}}
\newcommand{\Fqw}{{{\Bbb F}^\omega_q}}
\newcommand{\Md}{{{\Bbb M}_d}}
\newcommand{\Lt}{L\{\tau\}}
\newcommand{\Ct}{C\{\tau\}}
\newcommand{\Ctt}{C\{\tau,\tau^{-1}\}}
\newcommand{\mm}{{\frak m}}
\newcommand{\ttt}{{\frak t}}
\newcommand{\pp}{{\frak p}}
\newcommand{\qq}{{\frak q}}
\newcommand{\FF}{{\cal F}}
\newcommand{\GG}{{\cal G}}
\newcommand{\PP}{{\cal P}}
\newcommand{\Ffrakq}{{{\Bbb F}_\qq}}
\newcommand{\defequal}{\stackrel{\text{def}}{=}}
\newcommand{\nichts}{$\left.\right.$}

\begin{document}

\title[Fractional Power Series]{Fractional Power Series and Pairings on Drinfeld Modules}
\author{Bjorn Poonen}
\thanks{This research was supported by a Sloan Doctoral Dissertation Fellowship and an NSF Mathematical Sciences Postdoctoral Research Fellowship.  Research at MSRI is supported in part by NSF grant DMS-9022140.}
\address{Mathematical Sciences Research Institute \\ Berkeley, CA 94720-5070, USA}
\email{poonen@@msri.org}
\date{May 18, 1995}

\maketitle

\begin{abstract}
Let $C$ be an algebraically closed field containing $\Fq$ which is complete with respect to an absolute value $|\;|$.
We prove that under suitable constraints on the coefficients, the series $f(z) = \sum_{n \in \Z} a_n z^{q^n}$
converges to a surjective, open, continuous $\Fq$-linear homomorphism $C \rightarrow C$ whose kernel is locally compact.
We characterize the locally compact sub-$\Fq$-vector spaces $G$ of $C$ which occur as kernels of such series, and describe the extent to which $G$ determines the series.
We develop a theory of Newton polygons for these series which lets us compute the Haar measure of the set of zeros of $f$ of a given valuation, given the valuations of the coefficients.
The ``adjoint'' series $f^\ast(z) = \sum_{n \in \Z} a_n^{1/q^n} z^{1/q^n}$
converges everywhere if and only if $f$ does, and in this case there is a natural bilinear pairing
	$$\ker f \times \ker f^\ast \rightarrow \Fq$$
which exhibits $\ker f^\ast$ as the Pontryagin dual of $\ker f$.
Many of these results extend to non-linear fractional power series.
We apply these results to construct a Drinfeld module analogue of the Weil pairing, and to describe the topological module structure of the kernel of the adjoint exponential of a Drinfeld module.
\end{abstract}

\section{Introduction}
\label{intro}

Let $L$ be a perfect field of characteristic $p$.
Recall that a polynomial $f(x) \in L[x]$ is called {\em additive} if $f(x+y)=f(x)+f(y)$ identically.
It is easy to see that a polynomial is additive if and only if it is of the form
	$$f(x)=a_0 x + a_1 x^p + a_2 x^{p^2} + \cdots + a_n x^{p^n}.$$
The set of additive polynomials forms a noncommutative ring in which $(f \circ g)(x)=f(g(x))$.
This ring is generated by the scalar multiplications $x \mapsto ax$ for $a \in L$, and by $\tau(x)=x^p$.
In this ring we can write the $f$ above as $a_0 + a_1 \tau + a_2 \tau^2 + \cdots + a_n \tau^n$, and we will denote the ring by $\Lt$, using braces instead of brackets to remind ourselves that it is a {\em twisted} polynomial ring in $\tau$.
The indeterminate $\tau$ does not commute with elements of $L$ (which are acting as scalar multiplications); instead if $a \in L$ then $\tau a = a^p \tau$ as additive polynomials.

Following Ore~\cite{Or1}, we define the {\em adjoint} of $f$ to be the expression
	$$f^\ast(x) = a_0 x + a_1^{1/p} x^{1/p} + a_2^{1/p^2} x^{1/p^2} + \cdots + a_n^{1/p^n} x^{1/p^n}.$$
This defines a function on $L$ (and its algebraic closure) since $L$ is perfect.
Goss observed that the kernel of $f$ (the set of zeros of $f$ in a fixed algebraic closure of $L$) generates the same field extension as the kernel of $f^\ast$, by using a formula of Ore~\cite{Or1} which expresses the zeros of $f^\ast$ explicitly in terms of determinants involving the zeros of $f$.
One might ask whether there is a natural Galois-equivariant isomorphism between $\ker f$ and $\ker f^\ast$.
In general, the answer is no.
But, as was discovered independently by the author and N. Elkies~\cite{elkies}, there is something which is just as good:

\begin{unnumbered}
There exists a natural Galois-equivariant pairing
	$$\ker f \times \ker f^\ast \rightarrow \Fp.$$
\end{unnumbered}

One of the main results of this paper is the generalization of this theorem to the situation where $f$ is an additive {\em power series}.
Let $C$ be an algebraically closed field of characteristic $p$ which is complete with respect to a (non-archimedean) absolute value $|\;|$.
If $f \in C[[z]]$ is an additive power series, its adjoint is a ``linear fractional power series'' of the form
	$$g(z) = b_0 z + b_1 z^{1/p} + b_2 z^{1/p^2} + \cdots$$
It is easy to see that such a series converges for all $z$ when the $b_i$ tend to zero.
In this case, we will show that $g$ defines a continuous group homomorphism $C \rightarrow C$, which is surjective and open if nonzero.
One can see immediately that if $b_0 \not = 0$ then $g(z) = b_0 z + \text{(smaller terms)} \not= 0$ for large $z$, so $\ker g$ is a closed and bounded subgroup of $C$.
In fact, we prove that $\ker g$ is compact, and if infinitely many $b_i$ are nonzero, it is isomorphic to $\Fpw$ as a topological group.
Conversely, if $G$ is any compact subgroup of $C$, then $G$ is the kernel of such a $g$ if and only if the integral $\int_G v$ of the real-valued function $v(z)=-\log|z|$ with respect to Haar measure on $G$ equals $+\infty$.
Also, $G$ essentially determines $g$.
We next develop the theory of Newton polygons for such series.
Finally we prove the analogue of our pairing theorem above, namely that for any everywhere convergent additive power series $f \in C[[z]]$, there is a natural pairing
	$$\ker f \times \ker f^\ast \rightarrow \Fp$$
which exhibits the compact group $\ker f^\ast$ as the Pontryagin dual of the discrete group $\ker f$.

Actually we will generalize further by considering $\Fq$-linear power series for any power $q$ of $p$, and more substantially by considering {\em bi-infinite series}
	$$f(z) = \cdots + a_{-2} z^{1/q^2} + a_{-1} z^{1/q} + a_0 z + a_1 z^q + a_2 z^{q^2} + \cdots.$$
The kernels of these series are locally compact.
In fact, many of our results extend to non-linear series for which the set of allowable exponents is
	$$E=\{mp^n : m,n \in \Z, m \ge 0 \}.$$
We will develop results in this general context, specializing when appropriate.

Finally we give some applications of our results to the theory of Drinfeld modules.
We construct a pairing for Drinfeld modules which behaves in many ways like the Weil pairing on abelian varieties.
In addition, our results let us describe the topological module structure of the kernel of the adjoint exponential function of a Drinfeld module.
We conclude with a few unanswered questions.

\section{The rings $\PP$ and $\FF$ of fractional power series}
\label{tworings}

Let $C$ be an algebraically closed field containing $\Fq$ which is complete with respect to a (non-archimedean) absolute value $|\;|$.
Let $V = \{ x \in C : |x| \le 1 \}$ be the valuation ring of $C$, let $\mm = \{ x \in C : |x|<1 \}$ be its maximal ideal, and let $k=V/\mm$ be the residue field.

We will consider ``fractional power series'' in which the exponents belong to
	$$E=\{mp^n : m,n \in \Z, m \ge 0 \}.$$
Since we have unique $p$-th roots in $C$, for any $z \in C$ and $e = mp^n \in E$ we can interpret $z^e$ as $(z^{1/p^n})^m$.
Define $\PP$ to be the set of series $\sum_{e \in E} a_e z^e$ with coefficients $a_e$ in $C$, which converge at all $z \in C$.
Since $E$ is countable, it is clear what convergence means: for each $z \in C$ and $r>0$, there must be only finitely many terms in the series of absolute value greater than $r$.
And since we are in the non-archimedean situation, convergence at some $z$ of absolute value $r$ implies uniform convergence on all of $\Delta_r$.

Furthermore, we define the set $\FF$ of ``$\Fq$-linear fractional power series'' as the set of bi-infinite series
	$$f(z) = \cdots + a_{-2} z^{1/q^2} + a_{-1} z^{1/q} + a_0 z + a_1 z^q + a_2 z^{q^2} + \cdots$$
which converge for all $z \in C$.
Hence $\FF$ is a subset of $\PP$.
For $f \in \FF$ we will also write
	$$f = \cdots + a_{-2} \tau^{-2} + a_{-1} \tau^{-1} + a_0 + a_1 \tau + a_2 \tau^2 + \cdots,$$
thinking of $\tau$ as the operator $\tau(z)=z^q$ on $C$, and thinking of $a \in C$ acting as the scalar multiplication by $a$.
Note that for $a \in C$, $\tau a = a^q \tau$ and $\tau^{-1} a = a^{1/q} \tau^{-1}$, as maps on $C$.
We also define $\FF^+$ as the set of $f \in \FF$ such that $a_n=0$ for $n<0$ (these are just the convergent $\Fq$-linear power series), and $\FF^-$ as the set of $f \in \FF$ such that $a_n=0$ for $n>0$.

There is a simple criterion for the convergence of linear series:

\begin{prop}
\label{convergence}
Suppose $a_i \in C$ for $i \in \Z$.
The bi-infinite series
	$$f(z) = \cdots + a_{-2} z^{1/q^2} + a_{-1} z^{1/q} + a_0 z + a_1 z^q + a_2 z^{q^2} + \cdots$$
converges for all $z \in C$ if and only if the following two conditions hold:
\begin{enumerate}
	\item $\lim_{n \rightarrow -\infty} |a_n| = 0$.
	\item $\lim_{n \rightarrow +\infty} |a_n|^{1/q^n} = 0$.
\end{enumerate}
\end{prop}

\begin{pf}
By the well-known formula for the radius of convergence of a power series, the right end of the series converges everywhere when $\lim_{n \rightarrow +\infty} |a_n|^{1/q^n} = 0$.
The left end converges at $z$ when $\lim_{n \rightarrow -\infty} |a_n z^{q^n}| = 0$.
But for $z$ nonzero, this is equivalent to $\lim_{n \rightarrow -\infty} |a_n| = 0$, since $\lim_{n \rightarrow -\infty} |z^{q^n}| = 1$.
\end{pf}

For $f,g \in \PP$, we define $f+g$ and $f \cdot g$ in the obvious way, by expanding and grouping terms with the same exponent.

\begin{prop}
$(\PP,+,\cdot)$ is a commutative ring.
\end{prop}

\begin{pf}
The only nontrivial part is to check that the sum and product are actually everywhere convergent.
In fact, this is easy as well, since the sum or product of two convergent series is a doubly infinite convergent series, even before grouping terms.
\end{pf}

If $f,g \in \FF$, we can still use the definitions of addition and multiplication above, but $f \cdot g$ need not be an element of $\FF$, so $\FF$ is not a subring of $\PP$.
Nevertheless, we will make $\FF$ into a ring by using composition in place of multiplication.
In fact, we will define the composition $f \circ g$ of two general elements $f(z)=\sum_{d \in E} a_d z^d$ and $g(z)=\sum_{e \in E} b_e z^e$ of $\PP$.
First, for integral $m \ge 0$, let $g^m$ denote $g$ multiplied by itself $m$ times.
For $n \in \Z$, define $g^{p^n} = \sum_{e \in E} b_e^{p^n} z^{ep^n}$, which is consistent with the previous sentence.
For $d=mp^n \in E$, define $g^d = \left( g^m \right)^{p^n}.$
Finally define
	$$f \circ g = \sum_{d \in E} a_d g^d,$$
in which we expand and group terms of the same exponent.
If $f=\sum_{i \in \Z} a_i \tau^i$ and $g=\sum_{j \in \Z} b_j \tau^j$ are elements of $\FF$, then the above definition simplifies to $f \circ g =\sum_{n \in \Z} c_n \tau^n$ where
	$$c_n = \sum_{i+j=n} a_i b_j^{q^i}.$$

\begin{prop}
\label{composition}
For $f,g \in \PP$, $f \circ g$ converges to an element of $\PP$, and represents the composition of the maps $f$ and $g$.
\end{prop}

\begin{pf}
Fix $r>0$.
Since the series $g$ converges uniformly on $\Delta_r$, we may let $R = \sup |g(\Delta_r)|$.
Fix $z$ of absolute value $r$.
Then from the definition of $g^d$ it is clear that all terms in the expansion of $g^d$ (for fixed $d \in E$) are bounded in absolute value by $R^d$, with only finitely terms of absolute value greater than $s$, for each $s>0$.
Fix $t>0$.
There are only finitely many $d$ such that $|a_d|R^d > t$ (since $f$ is convergent), so there are only finitely many $d$ such that $a_d g^d$ has terms of absolute value greater than $t$, and the previous sentence implies that even for such $d$, there are at most finitely many such terms.
Hence even before grouping terms, the series for $f \circ g$ converges at $z$, and this implies that the grouping gives a convergent series for each coefficient of $f \circ g$ as well.
\end{pf}

\begin{prop}
\label{ring}
$(\FF,+,\circ)$ is a noncommutative ring containing $C$.
The center of $\FF$ is $\Fq$.
\end{prop}

\begin{pf}
By the previous proposition and the explicit form of the definition for linear series, if $f,g \in \FF$, then $f \circ g \in \FF$.
Hence the first part is clear.
The center of $\FF$ clearly contains $\Fq$.
On the other hand, if $c \in C$ is transcendental, direct computation shows that its centralizer in $\FF$ is $C$, and the only elements of $C$ which commute with $\tau$ are those in $\Fq$.
\end{pf}

\begin{prop}
\label{uniformcontinuity}
Any $f \in \PP$ defines a continuous function $f: C \rightarrow C$.
The continuity is uniform on each bounded subset of $C$.
\end{prop}

\begin{pf}
It suffices to show $f$ is uniformly continuous on $\Delta_R$ for each $R>0$.
Given $\epsilon >0$, there are only finitely many terms $a_e z^e$ in $f$ for which $|a_e| R^e >\epsilon$, since $f$ is everywhere convergent.
The sum of these is clearly uniformly convergent on $\Delta_R$, and the sum of the rest of the terms is bounded in absolute value by $\epsilon$ for all $z \in \Delta_R$, so $f$ is uniformly continuous on $\Delta_R$ as well.
\end{pf}

\begin{prop}
\label{nonvanishing}
If $f \in \PP$ is nonzero as a series, then the function it defines on $C$ is not identically zero.
\end{prop}

\begin{pf}
Since $f(1)$ converges, there are only finitely many coefficients of $f$ larger than a given size, and by scaling by an element of $C$, we may assume the largest coefficients have absolute value 1.
Then the reduction of $f$ modulo $\mm$ is a fractional {\em polynomial} over the residue field $k$, i.e., the $p^n$-th root of a polynomial.
Since $k$ is infinite, this reduction is nonzero at some $\bar{\alpha} \in k$.
If $\alpha \in C$ is any lift of $\bar{\alpha}$, then $|f(\alpha)|=1$, so in particular $f(\alpha) \not= 0$.
\end{pf}

\begin{prop}
\label{additive}
If $f \in \PP$ considered as a map on $C$ is $\Fq$-linear, then $f \in \FF$.
\end{prop}

\begin{pf}
Write $f= \sum_{e \in E} a_e z^e$.
For each $c \in C$, the function $f(cz+z)-f(cz)-f(z)$ in $\PP$ is identically zero, so by Proposition~\ref{nonvanishing} it is zero as a series.
Hence for each $e=mp^n$, the coefficient of $z^e$ in it is zero.
This coefficient is $a_e \left[ (c+1)^e - c^e - 1 \right]$.
If $e$ is not a power of $p$ and $a_e \not=0$, then the $p^n$-th power of this coefficient will be a nonzero polynomial in $c$, and hence will not vanish for some value of $c$.
This contradiction shows that $a_e=0$ for $e$ not a power of $p$.
A similar argument with $f(cz)-cf(z)$ for $c \in \Fq$ shows that $a_e=0$ for $e$ not a power of $q$; i.e., $f \in \FF$.
\end{pf}

\section{The topology and norm on $\PP$ and $\FF$}
\label{topologynorm}

We give $\PP$ the ``bounded-open topology.''
For real $r>0$, define
\begin{eqnarray*}
	\Delta_r	& = &	\{\, x \in C : |x| \le r \,\},	\\
	\Delta_r'	& = &	\{\, x \in C : |x| < r \,\}.
\end{eqnarray*}
If $r,s>0$ let
	$$\PP(r,s)=\{\,f \in \PP | f(\Delta_r) \subset \Delta_s \,\}.$$
Take these subsets of $\PP$ as a subbasis for the neighborhoods of zero.
For each $\alpha,\beta \in C$ and $c \in \Fq$, the set of $f \in \PP$ such that $f(\alpha+\beta)=f(\alpha)+f(\beta)$ and the set of $f \in \PP$ such that $f(c \alpha) = c f(\alpha)$ are closed subsets of $\PP$.
Hence by Proposition~\ref{additive}, $\FF$ is a closed subset of $\PP$.
We give $\FF$ the subspace topology, which is also the bounded-open topology on $\FF$.

For $f = \sum_{e \in E} a_e z^e \in \PP$, define
	$$\|f\| = \sup\; \{\, |a_e|:e \le 1 \,\} \cup \{\, |a_e|^{1/e}:e \ge 1 \,\}.$$
The convergence of $f$ implies that this is a finite real number, and that the supremum is attained.
If $f = \sum_{n \in \Z} a_n \tau^n \in \FF$ this becomes
	$$\|f\| = \sup\; \{\, |a_n|:n \le 0 \,\} \cup \{\, |a_n|^{q^{-n}}:n \ge 0 \,\}.$$
This norm satisfies an ultrametric triangle inequality: $\|f+g\| \le \max\{\|f\|,\|g\|\}$ for any $f,g \in \PP$.
In showing that the norm is consistent with the bounded-open topology on $\PP$ we will use the following.

\begin{lemma}
\label{termwise}
Let $f = \sum_{e \in E} a_e z^e \in \PP$, and let $r,s>0$.
Then the following are equivalent:
\begin{enumerate}
	\item $f \in \PP(r,s)$.
	\item $a_e z^e \in \PP(r,s)$ for all $e \in E$.
	\item $|a_e| r^e \le s$ for all $e \in E$.
\end{enumerate}
\end{lemma}

\begin{pf}
It is clear that~(2) and~(3) are equivalent, and that they imply~(1).
In showing~(1) implies~(3), we can assume $r$ is the absolute value of some element of $C$, because if not, it is an increasing limit of such $r$, and the implication follows if it is known for those $r$.
Then we can reduce to the case $r=1$ by composing $f$ with a scalar multiplication on the right.
Now what we must show is that if $f$ is a nonzero element of $\PP(1,s)$, a largest coefficient $b$ of $f$ is at most $s$ in absolute value.
The mod $\mm$ reduction of $b^{-1} f$ is the $p^n$-th root of a nonzero {\em polynomial}, so we can pick $x \in V$ which does not reduce to one of its roots, and find that $|b^{-1} f(x)|=1$.
Hence $|b| = |f(x)| \le s$, since $f \in \PP(1,s)$.
\end{pf}

\begin{prop}
\label{agree}
The bounded-open topology on $\PP$ is the same as that induced by the norm $\|\;\|$.
\end{prop}

\begin{pf}
It is clear from the definitions that if $0<s<1<r$,
	$$\|f\| \le s/r \implies f \in \PP(r,s).$$
On the other hand, by Lemma~\ref{termwise}, if $0<s<1$,
	$$f \in \PP(1,s) \cap \PP(1/s,1) \implies \|f\| \le s.$$
\end{pf}

\begin{theorem}
\label{topologicalring}
$(\PP,+,\cdot)$ is a complete topological ring, and the composition map $\PP \times \PP \rightarrow \PP$ is continuous.
Also, $(\FF,+,\circ)$ is a complete topological ring.
\end{theorem}

\begin{pf}
The triangle inequality for $\|\;\|$ implies that addition is continuous.
To check that multiplcation is continuous, we must show that if $f,g \in \PP$, and $r,s>0$, then for $f_0,g_0$ sufficiently small,
	$$(f+f_0)\cdot(g+g_0)-f\cdot g \in \PP(r,s).$$
Since $f$ and $g$ are bounded on $\Delta_R$, this follows for $f_0,g_0 \in \PP(r,\epsilon)$ for sufficiently small $\epsilon>0$.

Now let us show that $\PP$ is complete.
It suffices to show that if $f_1,f_2,\ldots$ is a sequence in $\PP$ tending to zero, then $\sum_{i=1}^\infty f_i$ converges to another element of $\PP$.
For each $e \in E$, the coefficients of $z^e$ in $f_i$ tend to zero by Lemma~\ref{termwise}, so they sum to some $a_e \in C$.
Interchanging the order of summation shows that $\sum_{i=1}^\infty f_i$ converges to $f \defequal \sum_{e \in E} a_e z^e$ uniformly on bounded subsets, and that $f \in \PP$.

For the continuity of composition, we must show that for $f_0,g_0$ sufficiently small,
\begin{equation}
\label{compcont}
	(f+f_0) \circ (g+g_0) - f \circ g \in \PP(r,s).
\end{equation}
Let $R=\sup |g(\Delta_r)|$.
If $g_0 \in \PP(r,R)$ and $f_0 \in \PP(R,s)$, then
\begin{equation}
\label{firsthalf}
	f_0 \circ (g+g_0) \in \PP(r,s).
\end{equation}
By Proposition~\ref{uniformcontinuity}, $f$ is uniformly continuous on $\Delta_R$, so if $g_0 \in \PP(r,\epsilon)$ for sufficiently small $\epsilon$,
\begin{equation}
\label{secondhalf}
	f \circ (g+g_0) - f \circ g \in \PP(r,s).
\end{equation}
Adding~(\ref{firsthalf}) and~(\ref{secondhalf}) yields~(\ref{compcont}).

By restriction from $\PP$ to $\FF$, the ring operations on $\FF$ are continuous.
As remarked earlier, $\FF$ is a closed subset of $\PP$, so the completeness of $\FF$ follows from that of $\PP$.
\end{pf}

Although the ring operations for both $\PP$ and $\FF$ are continuous with respect the topology induced by $\|\;\|$, it is not true that $\|f \cdot g\| \le \|f\| \cdot \|g\|$ for all $f,g \in \PP$; nor is it true that $\|f \circ g\| \le \|f\| \cdot \|g\|$ for all $f,g \in \FF$.
In fact, these can fail even if $f \in C$.
The following proposition shows that this defect of $\|\;\|$ is unavoidable.
(This is different from the theory of Banach algebras, in which if one has a norm for which multiplication is continuous in each variable, one can define a new norm with the sub-multiplicative property.  See Section~3.1 of~\cite{kadison}.)

\begin{prop}
There is no norm $\|\;\|'$ on $\PP$ inducing the same topology as above and satisfying
\begin{enumerate}
	\item $\|f\|' \ge 0$ with equality if and only if $f=0$.
	\item $\|f \pm g\|' \le \|f\|' + \|g\|'$
	\item $\|1\|' = 1$.
	\item $\|f \cdot g\|' \le \|f\|' \cdot \|g\|'$.
\end{enumerate}
Similarly, there is no norm $\|\;\|'$ on $\FF$ inducing the same topology and satisfying (1), (2), (3), and (4) with $\cdot$ replaced by $\circ$.
\end{prop}

\begin{pf}
Suppose there were such a norm on $\PP$.
As $\epsilon \rightarrow 0$ in $C$, $\epsilon z \rightarrow 0$ in $\PP$, so $\|\epsilon z\|'<1$ for some nonzero $\epsilon \in C$.
Then the series $1 + \epsilon z + (\epsilon z)^2 + \cdots$ should converge in $\PP$ to an inverse of $1-\epsilon z$.
This is impossible since $1-\epsilon z$ is zero at $z=1/\epsilon$.

Similarly, since $1-\epsilon \tau$ has no inverse (with respect to composition), there cannot be a sub-multiplicative norm on $\FF$.
\end{pf}

\section{A Class of Locally Compact Subspaces of $C$}
\label{subspaces}

For $z \in C^\ast$, let
	$$v(z)=-\log |z| \in \R$$
be the (additive) valuation associated with $|\;|$.
By convention, $v(0)=+\infty$.

If $f \in \FF^+$, then as is well-known, $\ker f$ is a discrete sub-$\Fq$-vector space of $C$, which is ``concentrated near infinity'' in the sense that if it is infinite, it is countable and its elements tend to infinity.
If $f \in \FF^-$, it will turn out that $G \defequal \ker f$ is a {\em compact} sub-$\Fq$-vector space of $C$, this time ``concentrated near zero'' in the sense that the integral $\int_G v$ of the valuation with respect to Haar measure on $G$ is $+\infty$.
The kernel of a general (i.e., possibly infinite in both directions) element of $\FF$ will belong to a certain hybrid of these two classes of vector spaces, which we now define.

\begin{defn}
Fix a real number $r>0$.
Let $\GG$ be the collection of sub-$\Fq$-vector spaces $G$ of $C$ satisfying the following conditions:
\begin{enumerate}
	\item $G \cap \Delta_r$ is compact.
	\item $\int_{G \cap \Delta_r} v = +\infty$.  The integral is with respect to (any) Haar measure on $G \cap \Delta_r$.
	\item $G \cap \Delta_r$ has finite or countable index in $G$, and in the latter case, if $g_1,g_2,\ldots$ is a system of coset representatives, then $\lim_{n \rightarrow \infty} |g_n| = +\infty$.
\end{enumerate}
\end{defn}

\begin{defn}
Let $G$ be a sub-$\Fq$-vector space of $C$.
Let $\{\lambda_n\}_{n \in \Z}$ be a bi-infinite collection of elements of $C$, possibly terminating on either side.
We say $\{\lambda_n\}_{n \in \Z}$ is a {\em descending basis} for $G$ if
\begin{enumerate}
	\item The $\lambda_n$ are independent over $\Fq$ and the closure of their span is $G$.
	\item For each $m \in \Z$ (for which $\lambda_m$ exists), $\lambda_m$ is a smallest nonzero element of the $\Fq$-vector space spanned by $\{\, \lambda_n \mid n \le m \,\}$.
	\item The sequence terminates on the left, or $\lim_{n \rightarrow -\infty} |\lambda_n| = \infty$.
	\item The sequence terminates on the right, or $\sum_{n=1}^\infty v(\lambda_n) q^{-n} = +\infty$.
\end{enumerate}
\end{defn}

\begin{prop}
\label{class}
Let $G$ be a sub-$\Fq$-vector space of $C$.
Then
\begin{enumerate}
	\item The definition of $\GG$ does not depend on the choice of $r>0$.
	\item $G$ has a descending basis if and only if $G \in \GG$.
	\item If $G \in \GG$, then $G$ is locally compact.
	\item For $G$ discrete, $G \in \GG$ if and only if it is finite or countable, with its elements tending to infinity in the latter case.
	\item For $G$ compact, $G \in \GG$ if and only if $\int_G v = +\infty$.
	\item If $G \in \GG$, then any closed sub-$\Fq$-vector space $H$ of $G$ also belongs to $\GG$.
\end{enumerate}
\end{prop}

\begin{pf} \nichts

\noindent {\em Proof of (2): }
Suppose $G$ has a descending basis $\{\lambda_n\}$.
Then the $\lambda_n$ are decreasing in size by condition~(2) in the definition, and $G \cap \Delta_r$ is the closure of the span of the $\{\lambda_n : n \ge n_0\}$ where $\lambda_{n_0}$ is the first $\lambda$ of absolute value at most $r$.
If the sequence of $\lambda_n$ after $\lambda_{n_0}$ terminates, then $G \cap \Delta_r$ is finite, hence compact.
Otherwise, these $\lambda_n$ form a sequence tending to zero by condition~(4), and we get a topological isomorphism $\Fqw \rightarrow G \cap \Delta_r$ which sends $(a_0,a_1,a_2,\ldots)$ to $\sum_{i=0}^\infty a_i \lambda_{n_0+i}$.
(This, together with Theorem~\ref{locallycompact}, proves the remark in the introduction regarding the topological group structure of the kernel of an element of $\FF^-$.)
Hence $G \cap \Delta_r$ is compact in any case.
Thus condition~(1) of the definition of $\GG$ is verified for $G$.

Assume that the Haar measure on $G$ is normalized so that the open subspace $G \cap \Delta_r$ has measure 1.
If $G \cap \Delta_r$ is finite, then $\int_{G \cap \Delta_r} v = +\infty$ since $0$ has positive measure.
Otherwise, under the isomorphism above, the set of $(a_0,a_1,a_2,\ldots) \in \Fqw$ such that the first nonzero $a$ is $a_i$ corresponds to a subset of $G \cap \Delta_r$ of measure $q^{-i}-q^{-i-1}$, and the function $v$ takes the value $v(\lambda_{n_0+i})$ everywhere on this subset, since otherwise some nontrivial combination of the $\lambda$'s of the same size as $\lambda_{n_0+i}$ would be closer to zero than $\lambda_{n_0+i}$, violating condition~(2) in the definition of decreasing basis.
These sets cover $\Fqw$ (except for the point 0), so
	$$\int_{G \cap \Delta_r} v = \sum_{i=0}^\infty v(\lambda_{n_0+i}) (q^{-i}-q^{-i-1}).$$
This is $+\infty$, since it agrees up to a finite number of terms with
	$$(1-q^{-1})q^{n_0} \sum_{n=1}^\infty v(\lambda_n) q^{-n}.$$
(Substitute $n=n_0+i$.)
Hence condition~(2) of the definition of $\GG$ holds.

Also, the finite linear combinations of the $\lambda_n$ for $n<n_0$ are coset representatives for $G \cap \Delta_r$ in $G$, so condition~(3) of the definition of $\GG$ follows from properties (2) and~(3) of a decreasing basis.
Hence $G \in \GG$.

\medskip
Conversely, suppose $G \in \GG$.
Since $G \cap \Delta_r$ is compact, only zero can be a limit point of $\{|x|:x \in G \cap \Delta_r\}$.
Combined with condition~(3) of the definition of $\GG$, this implies that the nonzero absolute values of elements of $G$ consist of a decreasing sequence $\{r_n\}_{n \in \Z}$, possibly terminating in either direction.
For each $r_n$, choose an $\Fq$-basis for $(G \cap \Delta_{r_n})/(G \cap \Delta_{r_n}')$, which is finite, by condition~(1) in the definition of $\GG$ when $r_n \le r$, and by condition~(3) when $r_n>r$.
We claim that the concatenation of these bases (lifted to $C$) is a descending basis for $G$.
Properties (1) and~(2) of a descending basis are clearly satisfied.
Property~(3) follows from condition~(3) in the definition of $\GG$, and property~(4) is equivalent to condition~(2) in the definition of $\GG$, as in the proof of the converse above.

\medskip
\noindent {\em Proof of (1), (3), (4), and (5): }
Part~(1) follows trivially from~(2).
Part~(3) follows since if $G \in \GG$, the open subgroup $G \cap \Delta_r$ is compact.
For part~(4), we may assume by~(1) that $r$ is smaller than the absolute value of the smallest nonzero element of $G$, and see what the conditions in the definition of $\GG$ say in this case.
Similarly, part~(5) follows by assuming $r$ is greater than the absolute value of all elements of $G$.

\medskip
\noindent {\em Proof of (6): }
In constructing a descending basis for $G$ as in the proof of~(2), we may assume that the chosen basis for $(G \cap \Delta_{r_n})/(G \cap \Delta_{r_n}')$ contains a basis for its subspace $(H \cap \Delta_{r_n})/(H \cap \Delta_{r_n}')$.
It is then easy to check that the concatenation of the latter bases (lifted to $C$) is a descending basis for $H$, so that $H \in \GG$ by~(2).
\end{pf}

\section{Maps and Fibers}
\label{kernels}

Here we state basic results on the maps on $C$ defined by elements of $\PP$ or $\FF$, and in particular study their fibers.
The proofs will be postponed until the next section.
First we have a fractional version of ``Picard's Theorem,'' which in the non-archimedean situation asserts that non-constant everywhere convergent power series are surjective.
The openness of the map is a kind of local surjectivity.

\begin{theorem}
\label{nonlinearsurjective}
If $f \in \PP$ is not a constant series, then $f:C \rightarrow C$ is a surjective open map.
\end{theorem}

This theorem has two easy consequences.

\begin{cor}
\label{ffdomain}
If $f,g \in \PP$ are not constant, then $f \circ g$ is not constant.
In particular, the noncommutative ring $\FF$ has no zero divisors.
\end{cor}

\begin{cor}
\label{distinct}
Distinct series of $\PP$ yield distinct maps on $C$.
\end{cor}

\begin{theorem}
\label{locallycompact}
If $f = \sum_{e \in E} a_e z^e \in \PP$ is not constant, then for each $c \in C$, $f^{-1}(c)$ is a locally compact subset of $C$, and is compact if and only if there is a largest $e$ for which $a_e \not=0$.
If $f = \sum_{n \in \Z} a_n \tau^n \in \FF$ is nonzero, then $G \defequal \ker f$ is in $\GG$.
Moreover,
\begin{center}
	$G$ is discrete $\iff$ $a_n=0$ for $n \ll 0$.	\\
	$G$ is compact $\iff$ $a_n=0$ for $n \gg 0$.
\end{center}
\end{theorem}

\noindent (The last two statements are dual in the sense of Section~\ref{adjointpd}.)
We have a converse for $f \in \FF$:

\begin{theorem}
\label{existence}
Each $G \in \GG$ is the kernel of some $f \in \FF$.
The $f$ is essentially unique: any $g \in \FF$ with the same kernel is of the form $\epsilon \tau^n \circ f$ for some $\epsilon \in C^\ast$ and $n \in \Z$.
(In other words, by Corollary~\ref{units} below, $g$ differs from $f$ only by a unit of $\FF$.)
\end{theorem}

Because $C$ is not locally compact, Theorem~\ref{locallycompact} implies it is not the union of $f^{-1}(0)$ and $g^{-1}(0)$ for nonzero $f,g \in \PP$.
This proves the following.

\begin{cor}
\label{ppdomain}
The commutative ring $\PP$ has no zero divisors.
\end{cor}

Recall that for a power series $f(z) = \sum_{i=0}^\infty a_i z^i \in C[[z]]$, the Newton polygon is the lower convex hull of the set of points $(i,v(a_i))$ in the plane, and that the Newton polygon gives information on the valuations of the zeros of $f$.
(See~\cite{amice},~\cite{artin},~\cite{kneser},~\cite{koblitz}, or~\cite{neukirch}.)
Similarly, we define the Newton polygon of $f=\sum_{e \in E} a_e z^e \in \PP$ to be the lower convex hull of the set of points $(e,v(a_e))$ in the plane.
Let $v_p$ denote the $p$-adic valuation on $\Q$.
For $f = \sum_{e \in E} a_e z^e \in \PP$, and $n \ge 0$, define
\begin{equation}
\label{fndefinition}
	f_n = \sum_{v_p(e) \ge -n} a_e z^e.
\end{equation}
Therefore $f_n$ is the $p^n$-th root of a power series.
Since $f$ is convergent everywhere, we see that the $f_n$ converge to $f$ uniformly on bounded subsets (hence also in $\PP$) as $n$ tends to infinity.

\begin{theorem}
\label{nonlinearnewton}
Suppose $f \in \PP$ is not constant, and $c \in C$.
Then there is a canonical measure $\mu$ on $f^{-1}(c)$, characterized by
	$$\mu \left( f^{-1}(c) \cap (\alpha+\Delta_r) \right) = \lim_{n \rightarrow \infty} \frac{\text{\# zeros of $(f_n-c)^{p^n}$ in $(\alpha+\Delta_r)$}}{p^n},$$
for all $\alpha \in C$, $r>0$.
(Here zeros are to be counted with multiplicity.)
The horizontal length of the segment of the Newton polygon of $f-c$ having slope $s$ (if any) equals
	$$\mu(\{\,z \in f^{-1}(c) \mid v(z)=-s \,\}).$$
Also,
	$$\mu(f^{-1}(c)) = \sup \{\,e \mid a_e \not=0 \,\}.$$
(Both might be $+\infty$.)
If $f \in C[[z]]$, then $\mu$ is the ``counting measure''; i.e., $\mu(\{x\})$ is the multiplicity of $x$ as a zero of $f-c$.
Finally, if $f \in \FF$ is nonzero and $c=0$, then $\mu$ is a Haar measure on the locally compact group $\ker f$.
\end{theorem}

\begin{cor}
\label{trivialkernel}
If $f \in \PP$ vanishes only at 0, then $f(z)=a z^e$ where $a \in C^\ast$ and $e \in E$.
\end{cor}

\begin{pf}
By Theorem~\ref{nonlinearnewton}, any series with at least two monomials has a ero of nonzero valuation.
\end{pf}

\begin{cor}
\label{units}
The unit group of $\FF$ consists of the elements of the form $a \tau^n$ where $a \in C^\ast$ and $n \in \Z$.
\end{cor}

\begin{pf}
Any unit of $\FF$ must have trivial kernel, and so must be of the form in Corollary~\ref{trivialkernel}, with the exponent a power of $q$.
On the other hand, it is easy to write down inverses of elements of this form.
\end{pf}

\begin{theorem}
\label{factors}
Let $f,g$ be nonzero elements of $\FF$.
Then $\ker g \subseteq \ker f$ if and only if there exists $h \in \FF$ such that $f=h \circ g$.
In this case, the $h$ is unique, and $\ker h = g(\ker f)$.
\end{theorem}

\section{Proofs}
\label{proofs}

Here we provide proofs of the theorems stated in the previous section.

\begin{lemma}
\label{mean}
Suppose $f(z)=\sum_{v_p(e) \ge -n} a_e z^e \in \PP$, and $j>0$ is such that $a_j \not= 0$.
Then there exists $x \in C$ such that $f(x)=0$ and $|x| \le |f(0)/a_j|^{1/j}$.
\end{lemma}

\begin{pf}
This follows from looking at the Newton polygon for the (ordinary) power series $f(z)^{p^n}.$
\end{pf}

Let $f_n$ be as in~(\ref{fndefinition}).

\begin{lemma}
\label{cosets}
Suppose $f \in \PP$ is not a constant series.
Fix $r,R \in \R$, with $0<r<R$, and fix $c \in C$.
Then the images of $f^{-1}(c) \cap \Delta_R$ and $f_n^{-1}(c) \cap \Delta_R$ in $\Delta_R/\Delta_r$ coincide for $n \gg 0$.
\end{lemma}

\begin{pf}
By considering $f-c$, we may assume $c=0$.
By composing with scalar multiplcations on both sides, we may assume $R=1$, and that the largest coefficient of $f$ other than the constant coefficient is of absolute value 1.
Let $a_e z^e$ be the term of $f$ with largest $e$ for which $|a_e|=1$.
Let $\epsilon_n=\sup|(f_{n+1}-f_n)(\Delta_R)|$.
Since $f_n \rightarrow f$, $\epsilon_n \rightarrow 0$.

Suppose $\lambda_n \in f_n^{-1}(0) \cap \Delta_R$.
Then
	$$|f_{n+1}(\lambda_n)| = |(f_{n+1}-f_n)(\lambda_n)| \le \epsilon_n.$$
Provided that $n+1 \ge -v_p(e)$, the coefficient of $z^e$ in $f_{n+1}(z+\lambda_n)$ has absolute value 1, by choice of $e$ and since $|\lambda_n| \le 1$.
Hence by Lemma~\ref{mean}, there exists $z_n$ such that $|z_n| \le \epsilon_n^{1/e}$ and
	$$f_{n+1}(z_n+\lambda_n)=0.$$
Let $\lambda_{n+1}=z_n+\lambda_n$.
Provided $n$ was chosen large enough, we can repeat the argument to construct a sequence $\lambda_n,\lambda_{n+1},\lambda_{n+2},\ldots$ with $\lambda_i \in f_i^{-1}(0)$ and
	$$|\lambda_{i+1}-\lambda_i| \le \epsilon_i^{1/e} \le r.$$
Since $\epsilon_i \rightarrow 0$, the $\lambda_i$ converge to a limit $\lambda_\infty$ with $|\lambda_\infty-\lambda_n| \le r$.
Also since $f_n \rightarrow f$ uniformly on $\Delta_R$,
	$$f(\lambda_\infty)=\lim_{i \rightarrow \infty} f_i(\lambda_i)=0.$$

For the other direction, suppose $\lambda \in f^{-1}(0) \cap \Delta_R$.
By Lemma~\ref{mean}, for $n \gg 0$, there exists $z_n$ such that $f_n(z_n+\lambda)=0$ and
	$$|z_n| \le |f_n(\lambda)|^{1/e} \rightarrow 0$$
as $n \rightarrow \infty$.
In particular, for sufficiently large $n$ (and how large does not depend on $\lambda$), $|z_n| \le r$, so $z_n+\lambda \in f_n^{-1}(0)$ is in the same coset of $\Delta_r$ as $\lambda$.
\end{pf}

From the preceding we can deduce two useful corollaries.
First of all, we see that Lemma~\ref{mean} extends to arbitrary non-constant $f \in \PP$.

\begin{cor}
\label{supermean}
Suppose $f(z)=\sum_{e \in E} a_e z^e \in \PP$, and $j>0$ is such that $a_j \not= 0$.
Then there exists $x \in C$ such that $f(x)=0$ and $|x| \le |f(0)/a_j|^{1/j}$.
\end{cor}

Next we have a finiteness result, which will be used to prove that fibers are locally compact.

\begin{cor}
\label{preimage}
If $f \in \PP$ is non-constant, $0<r<R$, and $c \in C$, then the image of $f^{-1}(c) \cap \Delta_R$ in $\Delta_R/\Delta_r$ is finite.
\end{cor}

\begin{pf}
The result with $f$ replaced by $f_n$ is proved by applying the theory of Newton polygons to the power series $f_n^{p^n}-c^{p^n}$.
Now apply Lemma~\ref{cosets}.
\end{pf}

\begin{pf*}{Proof of Theorem~\ref{nonlinearsurjective}}
Applying Corollary~\ref{supermean} to $f-c$ where $f \in \PP$ is not a constant series and $c \in C$ shows that there is a solution to $f(z)-c=0$, so $f$ is surjective.

In checking that $f$ defines an open map, we may reduce to proving $f$ is open at 0, by composing $f$ with a translation.
Moreover we may assume $f(0)=0$.
If $c \in C$ is small, Corollary~\ref{supermean} applied to $f-c$ shows that $f(z)-c=0$ has a solution near 0.
This is exactly what is needed to prove that $f$ is open at 0.
\end{pf*}

\begin{pf*}{Proof of Theorems \ref{locallycompact} and~\ref{nonlinearnewton}}
First let us show that for any $r>0$, $f^{-1}(c) \cap \Delta_r$ is compact.
Suppose $\alpha_1,\alpha_2,\ldots$ is a sequence in $f^{-1}(c) \cap \Delta_r$.
By Corollary~\ref{preimage}, these elements lie in finitely many cosets of $\Delta_{r/2}$, so we can find an infinite subsequence within one coset.
Next we can find a subsequence of this subsequence lying within a single coset of $\Delta_{r/3}$, and so on, with the $n$-th subsequence lying within a coset of $\Delta_{r/n}$.
By diagonalization we obtain a convergent subsequence of the original sequence.
Since $f$ is continuous, $f^{-1}(c)$ is closed (and so is $\Delta_r$), so the limit lies in $f^{-1}(c) \cap \Delta_r$.
Thus $f^{-1}(c) \cap \Delta_r$ is compact.
Taking larger and larger $r$ shows that $f^{-1}(c)$ is locally compact.

\medskip

Before completing the proof of Theorem~\ref{locallycompact}, let us turn to Theorem~\ref{nonlinearnewton}.
First we check that the limit in the definition of $\mu$ exists.
Without loss of generality we may assume $c=0$.
Also we may assume $\alpha=0$, by considering $f(z+\alpha)$.
Let $\ell_r(f)$ denote the total horizontal length of the segments of the Newton polygon of $f$ whose slope is less than $\log r$, which is the largest $x$-coordinate of a point of contact of the Newton polygon with a supporting line of slope $\log r$.
The fact that $f$ is everywhere convergent implies that below any line in $\R^2$ there are at most finitely vertices, so for each $r$, $\ell_r(f_n)=\ell_r(f)$ for $n \gg 0$.

The Newton polygon for $f_n^{p^n}$ is the dilation of that of $f_n$ by the factor $p^n$, so $\ell_r(f_n^{p^n})=p^n \ell_r(f_n)$.
The theory of Newton polygons for ordinary power series implies that $\ell_r(f_n^{p^n})$ counts the number of zeros of $f_n^{p^n}$ (with multiplicity) in $\Delta_r$, so the limit in the definition of $\mu$ (for $c=\alpha=0$) converges to $\ell_r(f)$.

To check that $\mu$ is truly a measure, it will suffice to check that the definition is consistent in the sense that when $r<R$, its value on $f^{-1}(c) \cap (\beta+\Delta_R)$ equals the sum of its values on $f^{-1}(c) \cap (\alpha+\Delta_r)$ with $\alpha+\Delta_r$ ranging over the cosets of $\Delta_r$ contained in $\beta+\Delta_R$.
By Lemma~\ref{cosets} and Corollary~\ref{preimage}, we need concern ourselves with only finitely many of these cosets.
Now the result is clear from the finite additivity of the right hand side of the definition of $\mu$.

By what we have shown so far,
	$$\mu(\{\, z \in f^{-1}(0) : |z| \le r \,\}) = \ell_r(f).$$
Taking this and subtracting the same with $r$ replaced by $r-\epsilon$ (with $s=\log r$, and $\epsilon$ sufficiently small) shows that
	$$\mu(\{\,z \in f^{-1}(c) \mid v(z)=-s \,\})$$
equals the horizontal length of the Newton polygon segment of slope $s$.

We now check the final statements of Theorem~\ref{nonlinearnewton}.
By what we have shown so far, $\mu(f^{-1}(c))$ is the sum of the horizontal lengths of all the segments, and this is clearly $\sup \{\,e \mid a_e \not=0 \,\}$.
If $f \in C[[z]]$, the theory of Newton polygons for power series tells us that the lengths of the horizontal segments are counting zeros, so $\mu$ is the counting measure.
Finally, if $f \in \FF$ is nonzero and $c=0$, it is clear from the definition that $\mu$ is translation-invariant on the locally compact group $\ker f$.
Moreover, $\mu(\ker f)>0$ by the formula we just derived for the measure of the whole space.
Thus $\mu$ is a Haar measure on $\ker f$.
This completes the proof of Theorem~\ref{nonlinearnewton}.

\medskip

We now resume the proof of Theorem~\ref{locallycompact}.
First, if there is a largest $e$ for which $a_e \not=0$, then the slopes of the Newton polygon of $f-c$ are bounded above for each $c \in C$, and hence by Theorem~\ref{nonlinearnewton}, $f^{-1}(c) \subset \Delta_r$ for some $r>0$.
Thus $f^{-1}(c)=f^{-1}(c) \cap \Delta_r$, which we showed already was compact.
On the other hand, if there is no largest $e$ for which $a_e \not=0$, then the Newton polygon of $f-c$ has infinitely many segments of increasing slope, so there is a sequence of elements in $f^{-1}(c)$ of decreasing valuation.
Such a sequence cannot have a convergent subsequence, so $f^{-1}(c)$ is not compact.

From now on, we assume $f=\sum_{n \in \Z} a_n \tau^n \in \FF$, $c=0$, and $G=\ker f$.
By the first paragraph of this proof, $G \cap \Delta_r$ is compact for any $r>0$, proving that $G$ satisfies condition~(1) in the definition of $\GG$.
Also by Corollary~\ref{preimage}, the image of $\ker f \cap \Delta_{nr}$ in $\Delta_{nr}/\Delta_r$ is finite for all $n$, so condition (3) in the definition of $\GG$ follows.

Now we check condition~(2); i.e., that if $r>0$, then $\int_{G \cap \Delta_r} v =+\infty$.
If we integrate only over the subset
	$$\{\, z \in G \mid v(z)=-s \,\}$$
we get $-s$ times the Haar measure of this set, which, by what we just proved, is the same as $-s$ times the horizontal length of the segment of the Newton polygon of slope $s$, which is the {\em vertical} displacement as one moves along the segment from right to left.
If we sum over all the segments of slopes less than $\log r$, we deduce that $\int_{G \cap \Delta_r} v$ equals the vertical displacement as one moves along the section of the Newton polygon to the left, starting at some point depending on $r$.
This vertical displacement is $+\infty$, since $v(a_n) \rightarrow +\infty$ as $n \rightarrow -\infty$.

Thus $G \in \GG$.
Finally, we check the last equivalences in Theorem~\ref{locallycompact}.
First, $G$ is discrete if and only if there are no zeros of large finite positive valuation.
By Theorem~\ref{nonlinearnewton}, this happens if and only if $a_n=0$ for $n \ll 0$.
Lastly, the criterion for $G$ to be compact follows from our earlier criterion for fibers of general elements of $\PP$.
\end{pf*}

For future reference, we record the following well-known result.

\begin{lemma}
\label{known}
If the zeros of a separable polynomial $f(z) \in C[z]$ form a sub-$\Fq$-vector space of $C$, then $f \in \Ct$.
\end{lemma}

\begin{pf}
Proposition~1.3 in~\cite{DH1} proves this for $q=p$, so $f$ is an additive polynomial.
Comparing zeros and linear coefficients shows that $f(cz)=cf(z)$ for all $z \in \Fq$, and this forces $f \in \Ct$.
\end{pf}

\begin{lemma}
\label{image}
If $G \in \GG$ and $f \in \FF$ has kernel $G \cap \Delta_r$ for some $r>0$,
then $f(G)$ is in $\GG$ and is discrete.
\end{lemma}

\begin{pf}
Since $f$ has compact kernel, its coefficient of $\tau^n$ is zero for $n \gg 0$ by Theorem~\ref{locallycompact}.
Hence in the series for $f(z)$, the last term dominates for large $z$.
Thus $|f(z)| \rightarrow \infty$ as $|z| \rightarrow \infty$.
Let $g_1,g_2,\ldots$ be representatives for the cosets of $G \cap \Delta_r$.
These form a sequence tending to infinity, by condition (3) in the definition of $\GG$.
Then
	$$f(G)=\{f(g_1),f(g_2),\ldots\}$$
also consists of a sequence tending to infinity, so it is in $\GG$ and is discrete.
\end{pf}

\begin{pf*}{Proof of Theorem~\ref{existence}, existence} \nichts

\noindent{\em Case 1:} $G$ discrete.

\noindent Then by Lemma~\ref{known}, we may take
	$$f(z) = z \prod_{g \in G, g \not=0} (1-z/g).$$

\noindent{\em Case 2:} $G$ compact.

\noindent We may assume $G$ is infinite, since otherwise we are in Case 1.
Without loss of generality suppose $G \subset \Delta_1$.
(If for some small $c \in C$ we can get $cG$ as a kernel of $f$, then we can get $G$ as the kernel of $f(cz)$.)

By Proposition~\ref{class}, there exists a descending basis $\lambda_1,\lambda_2,\ldots$ for $G$.
(It terminates on the left by condition~(3) of the definition of a descending basis.)
Let $V_n$ be the $\Fq$-vector space generated by $\lambda_1,\ldots,\lambda_n$.
Define
	$$g_n(z) = \prod_{\lambda \in V_n} (z-\lambda),
	\qquad \qquad h_n(z) = g_n(z)^{1/q^n}.$$
By Lemma~\ref{known}, $g_n \in \Ct$, so $h_n \in \FF^-$, and both have coefficients bounded in absolute value by 1, since $G \subset \Delta_1$.
We will eventually show that the $h_n$ converge to the desired $f \in \FF^-$.

We have
	$$g_{n+1}(z) = g_n(z)^q - g_n(\lambda_{n+1})^{q-1} g_n(z)$$
since both sides are monic $\Fq$-linear polynomials of degree $q^{n+1}$ which vanish on $V_{n+1}$.
Raising to the $1/q^{n+1}$ power yields 
\begin{equation}
\label{recurrence}
	h_{n+1}(z) = h_n(z) - c_n h_n(z)^{1/q}
\end{equation}
where $c_n=h_n(\lambda_{n+1})^{1-1/q}$.

We claim that $|c_n| \rightarrow 0$ as $n \rightarrow \infty$.
We have
\begin{eqnarray*}
	v(h_n(\lambda_{n+1}))	& = &	1/q^n \sum_{\lambda \in V_n} v(\lambda_{n+1}-\lambda)	\\
				& \ge &	1/q^n \sum_{\lambda \in V_n} \min\{v(\lambda_{n+1}),v(\lambda)\}	\\
				& \ge &	1/q^n \sum_{\lambda \in V_n, \lambda \not= 0} v(\lambda),
\end{eqnarray*}
since $v(\lambda_{n+1}) \ge v(\lambda)$ for all nonzero $\lambda \in V_n$, by the definition of descending basis.
As $n \rightarrow \infty$, this tends to $\int_G v = +\infty$.
Hence
	$$|c_n| = \exp(-(1-1/q) v(h_n(\lambda_{n+1}))) \rightarrow 0.$$

By~(\ref{recurrence}), we see that
	$$\|h_{n+1}(z) - h_n(z)\| = \|c_n h_n(z)^{1/q}\| \le |c_n| \rightarrow 0$$
as $n \rightarrow \infty$.
Hence the $h_n$ converge to some $f \in \FF$.

If $z \in V_n$, then $h_i(z)=0$ for $i \ge n$, so $f(z)=0$.
Thus $\ker f$ contains $\bigcup_{n=1}^\infty V_n$ and its closure, which is $G$.
On the other hand, if $z \not\in G$, and we set
	$$\delta = \inf_{\lambda \in G} |z-\lambda|,$$
then directly from the definitions of $g_n(z)$ and $h_n(z)$ we get
	$$|g_n(z)| \ge \delta^{q^n},	\qquad \qquad	|h_n(z)| \ge \delta,$$
for all $n$, so in particular, $f(z) = \lim_{n \rightarrow \infty} h_n(z) \not= 0$.
Hence $\ker f = G$, as desired.

\medskip
\noindent{\em Case 3:} $G \in \GG$ arbitrary.

\noindent By Case 2, we can find $g \in \FF$ with kernel $G \cap \Delta_1$.
By Lemma~\ref{image}, $g(G)$ is in $\GG$ and is discrete, so by Case 1, there exists $h \in \FF$ with kernel $g(G)$.
Let $f=h \circ g$.
Then
	$$f(x)=0 \iff h(g(x))=0 \iff g(x) \in g(G) \iff x \in G,$$
since $\ker g \subseteq G$.
\end{pf*}

Before proving Theorem~\ref{factors} and the uniqueness part of Theorem~\ref{existence}, let us prove the following ``remainder'' lemma.

\begin{lemma}
\label{remainder}
If $f \in \FF$, then there exists $q \in \FF$ such that
	$$f = q \circ (1 - \tau) + f(1).$$
\end{lemma}

\begin{pf}
Suppose $f = \sum_{n \in \Z} a_n \tau^n \in \FF$.
Without loss of generality, assume $f(1)=\sum_{n \in \Z} a_n =0$.
Then the rate of growth of the $a_n$ imposed by Proposition~\ref{convergence} implies that if
	$$b_n \defequal \sum_{m=-\infty}^n a_m = -\sum_{m=n+1}^\infty a_m,$$
then $q \defequal \sum_{n \in \Z} b_n \tau^n$ belongs to $\FF$.
(Use the first definition of $b_n$ to get convergence on the left, and the second for convergence on the right.)
Now
	$$q \circ (1-\tau) = \sum_{n \in \Z} (b_n-b_{n-1}) \tau^n = \sum_{n \in \Z} a_n \tau^n = f.$$
\end{pf}

\begin{pf*}{Proofs of Theorem~\ref{factors} and the uniqueness in Theorem~\ref{existence}}
In Theorem~\ref{factors}, it is clear that if $f=h \circ g$, then $\ker g \subseteq \ker f$, so we will concern ourselves with converse, that if $\ker g \subseteq \ker f$, then there exists $h \in \FF$ such that $f = h \circ g$.
Then the uniqueness of $h$ is clear from Corollary~\ref{ffdomain}, and $\ker h$ must be $g(\ker f)$ by the surjectivity of $g$ from Theorem~\ref{nonlinearsurjective}.

\noindent{\em Step 1: } Prove Theorem~\ref{factors} for $g= \alpha_0 \circ (1-\tau) \circ \alpha_1 \circ (1-\tau) \circ \alpha_2 \cdots (1-\tau) \circ \alpha_d$, where $\alpha_0,\alpha_1,\ldots,\alpha_d \in C^\ast$.

\medskip
\noindent We use induction on $d$.
The base case $d=0$ is trivial.
Suppose $d \ge 1$.
Let $\tilde{g} = \alpha_1 \circ (1-\tau) \circ \alpha_2 \cdots (1-\tau) \circ \alpha_d$.
Then the inductive hypothesis implies that $f = j \circ \tilde{g}$ for some $j \in \FF$.
Since $f$ kills $\ker g$, $j$ kills $\tilde{g}(\ker g) = \ker( \alpha_1 \circ (1-\tau)) = \Fq$.
In particular $j$ kills 1, so by Lemma~\ref{remainder}, $j=q\circ (1-\tau)$ for some $q$.
If we now let $h=q\circ \alpha_0^{-1}$, then $f=h \circ g$.

\medskip
\noindent{\em Step 2: } Every separable $g \in \Ct$ of $\tau$-degree $d$ is of the form in Step 1.

\medskip
\noindent Again we use induction on $d$.
The base case $d=0$ is trivial.
If $d \ge 1$, there exists $c \in \ker g$, $c \not=0$.
Then $g$ kills $\ker ((1-\tau) \circ c^{-1})$, so $g = h \circ (1-\tau) \circ c^{-1}$ for some $h \in \FF$ by Step 1, and the proof of the existence of $h$ in fact shows that $h$ must be a separable element of $\Ct$.
Applying the inductive hypothesis to $h$ produces the desired factorization.

\medskip
\noindent{\em Step 3: } Theorem~\ref{factors} holds when $g \in \Ctt$.

\medskip
\noindent Let $\tau^n$ be the lowest (most negative) power occuring in $g$.
Then by Steps 1 and 2, we can write $f=h' \circ (\tau^{-n} \circ g)$, so we can take $h=h' \circ \tau^{-n}$.

\medskip
\noindent{\em Step 4: } If $g_n,h_n \in \FF$ are nonzero, $f = h_n \circ g_n$ for each $n$, and the $g_n$ converge to some nonzero $g \in \FF$, then the $h_n$ converge to some $h \in \FF$ and $f = h \circ g$.

\medskip
\noindent Fix $x \in C$.
We have $h_n(x)=f(g_n^{-1}(x))$ where $g_n^{-1}(x)$ denotes any $y_n \in C$ such that $g_n(y_n)=x$.
(By Theorem~\ref{nonlinearsurjective}, such $y_n$ exist.)
We will construct a sequence of such $y_n$ which {\em converges}.
Fix a nonzero term $a_j \tau^j$ in $g$.
Choose some large $R>0$, then choose some large $n_0>0$.
(We'll specify how large as we go along; how large we need $n_0$ to be depends on how large $R$ was taken to be.)
If $n \gg n_0$, the coefficient of $\tau^j$ in $g_n$ has the same absolute value as $a_j$, provided $n_0$ was chosen large enough.
Then by Corollary~\ref{supermean}, there exists a solution $y_{n_0} \in \Delta_R$ to $g_{n_0}(y_{n_0})=x$, if $R$ was chosen large enough.
Now inductively define $\epsilon_n$ and $y_{n+1} \in \Delta_R$ for $n \ge n_0$ as follows.
Pick $\epsilon_n$ such that $g_{n+1}(\epsilon_n)=(g_n-g_{n+1})(y_n)$, and set $y_{n+1}=y_n+\epsilon_n$.
Since the $g_n$ converge uniformly on $\Delta_R$, and since $y_n \in \Delta_R$ by the inductive hypothesis, $(g_n-g_{n+1})(y_n)$ can be assumed to be arbitrary small, if $n_0$ was chosen large enough.
Then by Corollary~\ref{supermean}, $\epsilon_n \in \Delta_R$, $y_{n+1} \in \Delta_R$, and
\begin{eqnarray*}
	g_{n+1}(y_{n+1}) & = & g_{n+1}(y_n) + g_{n+1}(\epsilon_n) \\
			& = &  g_{n+1}(y_n) + g_n(y_n) - g_{n+1}(y_n) \\
			& = & x.
\end{eqnarray*}
Moreover Corollary~\ref{supermean} guarantees that the $\epsilon_n$ can be chosen going to zero, so the $y_n$ converge as desired.

Since $f$ is continuous, we see that the sequence $h_n(x)$ converges.
In fact, as is clear from the construction of the $y_n$, the rate of convergence depends not on $x$, but only on an upper bound for $|x|$, so the $h_n$ converge uniformly on every bounded subset of $C$.
By Theorem~\ref{topologicalring}, the $h_n$ converge to some $h \in \FF$.
Taking the limit of $f=h_n \circ g_n$ as $n \rightarrow \infty$ yields $f=h \circ g$.

\medskip
\noindent{\em Step 5: } Uniqueness holds in Theorem~\ref{existence} when $G \in \GG$ is compact or discrete.

Suppose $G \in \GG$ is compact.
By the proof of existence part of Theorem~\ref{existence}, there is a sequence of elements $j_n \in \Ctt$ with $\ker j_n \subseteq G$, which converges to an element $j \in \FF$ with kernel $G$.
Let $f$ be any other element of $\FF$ with kernel $G$.
By Step 3, we may write $f=h_n \circ j_n$, and then by Step 4, we may write $f=h \circ j$.
But $j$ is surjective by Theorem~\ref{nonlinearsurjective}, so $\ker h=j(\ker f)=0$ and hence $h$ is a unit by Corollaries \ref{trivialkernel} and~\ref{units}.

If $G \in \GG$ is discrete, we can find a sequence $j_n \in \Ct$ with $\ker j_n \subseteq G$, which converges to an element $j \in \FF$ with kernel $G$, so the same proof works.

\medskip
\noindent{\em Step 6: } Theorem~\ref{factors} holds when $\ker g$ is compact or discrete.

Indeed, by Step 5, $g$ must be a limit of elements $g_n \defequal h \circ j_n \in \Ctt$ with $\ker g_n = \ker j_n \subseteq \ker g$, and we can then apply Step 4.

\medskip
\noindent{\em Step 7: } Theorem~\ref{factors} holds.

By the existence part of Theorem~\ref{existence}, we can find $g_1 \in \FF$ with compact kernel $G \cap \Delta_1$.
By Step 6, we can write $g=j \circ g_1$, and $\ker j = g_1(\ker g)$, which is discrete by Lemma~\ref{image}.
By Step 6 again, we can write $f = h_1 \circ g_1$, and
	$$\ker j = g_1(\ker g) \subseteq g_1(\ker f) = \ker h_1.$$
By Step 6 yet a third time, we can write $h_1 = h \circ j$.
Then
	$$f = h_1 \circ g_1 = h \circ j \circ g_1 = h \circ g.$$

\medskip
\noindent{\em Step 8: } Uniqueness holds in Theorem~\ref{existence}.

Indeed, suppose $f$ and $g$ are nonzero elements of $\FF$ with the same kernel.
By Step 7, each of $f$ and $g$ is a left multiple of the other, so they differ by a unit.
\end{pf*}

\section{Pontryagin duality for $A$-modules}
\label{pda}

Here we recall and develop some results on locally compact topological modules to be used later in the paper.
Throughout, topological groups, rings, and modules are assumed Hausdorff, and each $\Hom$ consists of continuous homomorphisms and is given the compact-open topology.
If $G$ is a locally compact abelian topological group (LCA group), its {\em Pontryagin dual} is $\hat{G} \defequal \Hom_\Z(M,\R/\Z)$, which is again an LCA group.
The elements of $\hat{G}$ are called {\em characters} of $G$.
The main theorem of Pontryagin duality is that the natural map $G \rightarrow \Hat{\Hat{G}}$ is an isomorphism.
For an account of Pontryagin duality for groups, see~\cite{armacost}.

The following technical result on LCA groups is due to Kaplansky (unpublished) and first appears in a paper of Glicksberg~\cite{glicksberg}.
See Chapter~10 of~\cite{armacost} for a proof and discussion.

\begin{prop}
\label{topology}
Let $G$ be a group which becomes a LCA group under a topology $\ttt$.
Let $\ttt'$ be a strictly stronger locally compact group topology on $G$.
Then there is a $\ttt'$-continuous character $\gamma$ of $G$ which is not $\ttt$-continuous.
\end{prop}

\begin{cor}
If $f:G \rightarrow H$ is a surjective continuous homomorphism of LCA groups whose dual $\hat{f}:\hat{H} \rightarrow \hat{G}$ is also surjective, then $f$ is a topological isomorphism.
\end{cor}

\begin{pf}
Since $\hat{f}$ is surjective, $f$ is injective by~P.23(b) in~\cite{armacost}.
Thus $f$ is an isomorphism, except that $G$ may have a stronger topology than $H$.
If it were strictly stronger, then by the proposition above, there would be a character of $G$ not coming from a character of $H$, contradicting the surjectivity of $\hat{f}$.
Therefore $f$ is a topological isomorphism.
\end{pf}

\begin{cor}
\label{showdual}
Suppose $G$ and $H$ are LCA groups, and there is a continuous bilinear pairing
	$$G \times H \rightarrow \R/\Z$$
such that the induced maps
	$$f_1:G \rightarrow \hat{H} \;\;\;\;,\;\;\;\; f_2:H \rightarrow \hat{G}$$
are surjective.
Then $f_1$ and $f_2$ are topological isomorphisms.
\end{cor}

\begin{pf}
By a well-known property of the compact-open topology, the induced maps are continuous.
By Pontryagin duality (and chasing definitions), $f_2$ is the dual of $f_1$, so $f_1$ is a topological isomorphism by the previous corollary.
The symmetric argument works for $f_2$.
\end{pf}

We will make use of the theory of Pontryagin duality for topological modules developed by Flood~\cite{flood}.
Let $A$ be a locally compact commutative topological ring.
Consider the class $\A$ of locally compact topological $A$-modules.
If $M \in \A$, then the Pontryagin dual $\hat{M} \defequal \Hom_\Z(M,\R/\Z)$ of $M$ as a topological group has a natural $A$-module structure, so $\hat{M} \in \A$.
Let us restate the main theorem of~\cite{flood} for this situation:

\begin{theorem}[Pontraygin duality for $A$-modules]
\label{pdfora}
For each $M \in \A$, the topological $A$-modules $\hat{M}$ and $\Hom_A(M,\hat{A})$ are canonically isomorphic, and the canonical map
\begin{eqnarray*}
	M 	& \rightarrow	& \Hom_A(\Hom_A(M,\hat{A}),\hat{A})	\\
	m	& \mapsto	& (f \mapsto f(m))
\end{eqnarray*}
is an isomorphism of topological $A$-modules.
\end{theorem}

If $M,N \in \A$, we say $N$ is the Pontryagin dual of $M$ as an $A$-module, if it is isomorphic to $\Hom_A(M,\hat{A})$ as a topological $A$-module.
The following corollary is an $A$-module version of Corollary~\ref{showdual}.
The next corollary is useful in verifying that two elements of $\A$ are Pontryagin duals (as $A$-modules).

\begin{cor}
\label{showdual2}
Suppose $M,N \in \A$ and there is a continuous $A$-module pairing
	$M \times N \rightarrow \hat{A}$
such that the induced maps
	$$f_1:M \rightarrow \Hom_A(N,\hat{A}), \;\;\;\;\;\;\;\; f_2:N \rightarrow \Hom_A(M,\hat{A})$$
are surjective.
Then $f_1$ and $f_2$ are topological $A$-module isomorphisms.
\end{cor}

\begin{pf}
If we compose the pairing with the evaluation-at-1 map $\hat{A} \rightarrow \R/\Z$, we get a pairing
	$$M \times N \rightarrow \R/\Z$$
and the induced map $M \rightarrow \hat{N}$ is the same as $f_1$ once we identify $\hat{N}$ with $\Hom_A(N,\hat{A})$ as in Theorem~\ref{pdfora}.
Similarly $f_2$ is the same as the induced map $N \rightarrow \hat{M}$.
In particular, these induced maps are surjective, so by Corollary~\ref{showdual}, they are topological isomorphisms.  But we know that $f_1$ and $f_2$ are $A$-module homomorphism as well, so the result follows.
\end{pf}

\begin{lemma}
\label{showdual3}
Suppose $M,N$ are locally compact topological $\Fq$-vector spaces.
If there is a continuous $\Fq$-linear pairing
	$$\langle \;\; , \;\; \rangle : M \times N \rightarrow \Fq$$
such that the induced maps
	$$M \rightarrow \Hom_\Fq(N,\Fq) \;\;\;\;,\;\;\;\; N \rightarrow \Hom_\Fq(M,\Fq)$$
are surjective, then $N$ is the Pontryagin dual of $M$ as an $\Fq$-vector space.

If furthermore $A$ is a locally compact $\Fq$-algebra and $M$ and $N$ are topological $A$-modules such that $\langle am,n \rangle = \langle m,an \rangle$ for all $a \in A$, $m \in M$ and $n \in N$, then $\Hom_\Fq(A,\Fq) \cong \hat{A}$ as topological $A$-modules and the map
\begin{eqnarray*}
	[ \;\; , \;\; ] : M \times N	& \rightarrow &	\Hom_\Fq(A,\Fq) \cong \hat{A}	\\
	m \; , \; n \;\;		& \mapsto &	(a \mapsto \langle am,n \rangle)
\end{eqnarray*}
is a continuous $A$-module pairing which exhibits $N$ as the Pontryagin dual of $M$ as an $A$-module.
\end{lemma}

\begin{pf}
By counting elements, $\Fqhat$ is a one-dimensional $\Fq$-vector space, and we can explicitly identify $\Fq$ with $\Fqhat$ by mapping $1$ to the character $\chi=(1/p) \cdot \Tr_{\Fq/\Fp}$ of $\Fq$.
Thus the first half of the proposition is just the special case of Corollary~\ref{showdual2} for the discrete ring $\Fq$.

There is an isomorphism $\Hom_\Fq(A,\Fq) \cong \hat{A}$ of topological $\Fq$-vector spaces by Theorem~\ref{pdfora} (with $M=A$, $A=\Fq$), and it is easily checked that it preserves the $A$-module structure.

As in Corollary~\ref{showdual2}, the pairing $\langle \;\;,\;\; \rangle$ induces a pairing $M \times N \rightarrow \R/\Z$ which identifies $N$ with $\hat{M}$.
The condition $\langle am,n \rangle = \langle m,an \rangle$ ensures that this isomorphism is an isomorphism of $A$-modules.
Finally, if we identify $N \cong \hat{M}$ with $\Hom_A(M,\hat{A})$ using Theorem~\ref{pdfora}, we get a pairing
	$$M \times N \rightarrow \hat{A}$$
exhibiting $N$ as the Pontryagin dual of $M$ as an $A$-module, and definition chasing shows that this pairing is $[\;\;,\;\;]$.
\end{pf}

Now let $X$ be a nonsingular projective curve over $\Fq$, let $\infty$ be a closed point, and let $Y=X \setminus \infty$.
Let $A$ be the Dedekind ring of regular functions on the affine curve $Y$, and give $A$ the discrete topology.
Let $K$ be its fraction field, and let $K_\infty$ be the completion of $K$ at $\infty$.
The K\"{a}hler differentials of $A$ over $\Fq$ (see~\cite{matsumura}), i.e., the differentials on $X$ which are regular away from $\infty$, form a rank one projective $A$-module $\Omega_A$ inside the one-dimensional $K$-vector space $\Omega_K$ of K\"{a}hler differentials of $K$ over $\Fq$.
Therefore there is an ideal $J$ of $A$ isomorphic to $\Omega$ as an $A$-module.
Let $\Omega_\infty$ be the completion of $\Omega_K$ at $\infty$.
This is a one-dimensional vector space over $K_\infty$.
By Theorem~3 of Chapter~II in~\cite{weilnt}, any nontrivial character of the locally compact field $K_\infty$ can be used to identify $K_\infty$ with its dual, and it follows that the map
\begin{eqnarray}
\label{residuepairing}
	K_\infty \times \Omega_\infty	& \rightarrow &	\Fq	\\
	a \;\; , \;\; \omega \;\;	& \rightarrow &	\Res_\infty(a \omega) \nonumber
\end{eqnarray}
exhibits $K_\infty$ and $\Omega_\infty$ as Pontraygin duals of each other, as topological $\Fq$-vector spaces.

\begin{theorem}
\label{fractionaldual}
Let $I$ be a fractional ideal of $A$ (with the discrete topology).
Then the pairing~(\ref{residuepairing}) puts $I$ and $\Omega_\infty/(I^{-1} \Omega_A)$ in Pontryagin duality as $A$-modules, so the Pontryagin dual of $I$ as an $A$-module is $K_\infty/(I^{-1}J)$.
\end{theorem}

\begin{pf}
For $\omega \in \Omega_\infty$, we claim that $\omega \in \Omega_A$ if and only if $\Res_\infty(a \omega)=0$ for all $a \in A$.
For $n \ge 0$, let
\begin{eqnarray*}
	L(n \infty)	& = &	\{\, a \in K \mid \divv(a) \ge -n \infty \,\}	\\
	\Omega_{\ge n \infty}	& = &	\{\, \eta \in \Omega_\infty \mid \divv(\eta) \ge n \infty \,\}.
\end{eqnarray*}
(By the divisor of an element of $\Omega_\infty$ we mean only the part at $\infty$, which is the only part that makes sense.)
By the Theorem at the bottom of page~160 in~\cite{eichler}, with $h=0$, ${\frak a}=n \infty$, and $dw$ the principal part system which equals an element of $\Omega_K$ differing from $\omega$ by an element of $\Omega_{\ge n \infty}$ at $\infty$, and equalling zero at other places, we have that $\Res_\infty(a \omega)=0$ for all $a \in L(n \infty)$ if and only if $\omega \in \Omega_A + \Omega_{\ge n \infty}$.
But $A = \bigcup_{n \ge 0} L(n \infty)$, so $\Res_\infty(a \omega)=0$ for all $a \in A$ if and only if
	$$\omega \in \bigcap_{n \ge 0} (\Omega_A + \Omega_{\ge n \infty}) = \Omega_A,$$
since $\Omega_A$ is a discrete subgroup of $\Omega_\infty$, and the $\Omega_{\ge n \infty}$ form a decreasing neighborhood base of 0.

Hence $\Res_\infty(a \omega)=0$ for all $a \in I$ if and only if $\omega \in I^{-1} \Omega_A$.
By P.22(d) in~\cite{armacost}, the fact that the pairing~(\ref{residuepairing}) puts $K_\infty$ and $\Omega_\infty$ in Pontryagin duality implies that
	$$\hat{I} \cong \Hom_\Fq(I,\Fq) \cong \Omega_\infty/(I^{-1} \Omega_A) \cong K_\infty/(I^{-1}J),$$
and it is trivial to check that all our isomorphisms respect the $A$-module structures.
\end{pf}

Let $A_\qq$, $K_\qq$ denote the completions of $A$, $K$, respectively, at a nonzero prime $\qq$ of~$A$.

\begin{theorem}
\label{aqdual}
The $A_\qq$-module $\Hom_\Fq(K_\qq/A_\qq,\Fq)$ is free of rank one.
\end{theorem}

\begin{pf}
Let $t \in K$ be a uniformizing parameter at $\qq$.
Then $A_\qq \cong \Ffrakq[[t]]$, where $\Ffrakq$ is the residue field at $\qq$.
Now for the residue pairing
\begin{eqnarray*}
	\Ffrakq((t)) \times \Ffrakq((t))dt	& \rightarrow &	\Fq	\\
	a \;\;\; , \;\;\; \omega \;\;\;\;\;\;\;\;\;\;	& \mapsto &	\Res_\qq(a \omega)
\end{eqnarray*}
which exhibits the locally compact field $K_\qq$ as its own Pontryagin dual, the set of $\omega$ which give zero when paired with any $a \in \Ffrakq[[t]]$ is exactly $\Ffrakq[[t]]dt$, by definition of the residue, so by~P.22(d) in~\cite{armacost}, the Pontryagin dual of $A_\qq=\Ffrakq[[t]]$ is $K_\qq/A_\qq=\Ffrakq((t))/\Ffrakq[[t]]$.
\end{pf}

\section{Adjoints and Pontryagin duality}
\label{adjointpd}

The {\em adjoint map} on $\FF$ is the map $f \mapsto f^\ast$ defined by the following proposition.

\begin{prop}
There is a norm-preserving multiplication-reversing involution
\begin{eqnarray*}
		\FF	& \rightarrow	& \FF	\\
		f=\sum_{n \in \Z} a_n \tau^n	& \mapsto	& f^\ast=\sum_{n \in \Z} {a_n}^{1/q^n} \tau^{-n}.
\end{eqnarray*}
\end{prop}

\begin{pf}
The map is well-defined by Proposition~\ref{convergence}, and it is immediate from the definition that $\|f^\ast\|=\|f\|$.
One can check $f^{\ast \ast}=f$ and $(f \circ g)^\ast = g^\ast \circ f^\ast$ directly by comparing coefficients, or perhaps more enlighteningly by noting that
	$$f^\ast = \sum_{n \in \Z} \tau^{-n} a_n$$
is obtained from $f$ by interchanging $\tau$ and $\tau^{-1}$ and reversing the order of multiplication everywhere, and that this preserves the pair of commutation relations $\tau \circ a = a^q \circ \tau$ and $a \circ \tau^{-1} = \tau^{-1} \circ a^q$.
\end{pf}

We now construct an $\Fq$-linear pairing between $\ker f$ and $\ker f^\ast$, for nonzero $f \in \FF$.
Suppose $\alpha \in \ker f$ and $\beta \in \ker f^\ast$.
Then $f \circ \alpha$ vanishes on 1 (recall that $f \circ \alpha$ denotes the composition of $f$ with the map $x \mapsto \alpha x$), so by Lemma~\ref{remainder},
\begin{equation}
\label{def}
	f \circ \alpha = g_\alpha \circ (1 - \tau)
\end{equation}
for some $g_\alpha \in \FF$, which is uniquely defined, by Corollary~\ref{ffdomain}.
Define
	$$\langle \alpha,\beta \rangle_f = g_\alpha^\ast(\beta).$$

\begin{prop}
If $f \in \FF$ is nonzero, then
	$$\langle \;\;,\;\; \rangle_f : \ker f \times \ker f^\ast \rightarrow \Fq$$
is an $\Fq$-vector space pairing.
\end{prop}

\begin{pf}
First let us check that the image lands in $\Fq$.
Take adjoints of (\ref{def}):
\begin{equation}
\label{adjdef}
	\alpha \circ f^\ast = (1 - \tau^{-1}) \circ g_\alpha^\ast.
\end{equation}
Apply both sides to $\beta$:
	$$0 = (1 - \tau^{-1}) g_\alpha^\ast(\beta).$$
So
	$$g_\alpha^\ast(\beta) \in \ker(1 - \tau^{-1}) = \Fq.$$
From the definition~(\ref{def}), $g_{\alpha+\beta}=g_\alpha+g_\beta$, so the pairing is linear on the left.
Linearity on the right is obvious.
\end{pf}

We now prove a series of results leading up to Theorem~\ref{duality} below.

\begin{lemma}
\label{remainder2}
If $f \in \FF$, then there exists $q \in \FF$ such that
	$$f = (1 - \tau) \circ q + f^\ast(1).$$
\end{lemma}

\begin{pf}
Apply Lemma~\ref{remainder} to $f^\ast$ to get $h$ such that
	$$f^\ast = h \circ (1-\tau) + f^\ast(1).$$
Then taking adjoints yields
\begin{eqnarray*}
	f	& = &	(1-\tau^{-1}) \circ h^\ast + f^\ast(1)	\\
		& = &	(1-\tau) \circ (-\tau^{-1}) \circ h^\ast + f^\ast(1),
\end{eqnarray*}
so we may take $q=-\tau^{-1} \circ h^\ast$.
\end{pf}

\begin{lemma}
\label{oneone}
If $f,g \in \FF$ and $f \circ (1-\tau) = (1-\tau) \circ g$, then $f^\ast(1)=g(1)$.
\end{lemma}

\begin{pf}
Use Lemmas \ref{remainder} and~\ref{remainder2} to write
\begin{eqnarray*}
	f	& = &	(1-\tau) \circ q + f^\ast(1)	\\
	g	& = &	h \circ (1-\tau) + g(1).
\end{eqnarray*}
Substitute these into the given relation:
\begin{equation}
\label{subst}
	(1-\tau) \circ q \circ (1-\tau) + f^\ast(1) \circ (1-\tau) = (1-\tau) \circ h \circ (1-\tau) + (1-\tau) \circ g(1).
\end{equation}
Evaluating the given relation at 1 shows $g(1) \in \ker(1-\tau) = \Fq$, so $g(1)$ commutes with $1-\tau$ and we can cancel $1-\tau$ on the right in (\ref{subst}) to obtain
	$$(1-\tau) \circ q + f^\ast(1) = (1-\tau) \circ h + g(1).$$
Taking adjoints and evaluating at 1 gives the desired result.
\end{pf}

\begin{prop}
\label{backwards}
If $f \in \FF$ is nonzero, $\alpha \in \ker f$, and $\beta \in \ker f^\ast$, then $\langle \alpha,\beta \rangle_f = -\langle \beta,\alpha \rangle_{f^\ast}$.
\end{prop}

\begin{pf}
Write
\begin{eqnarray*}
	f \circ \alpha = g_\alpha \circ (1 - \tau)	\\
	f^\ast \circ \beta = h_\beta \circ (1 - \tau).
\end{eqnarray*}
Multiply the first by $\beta$ on the left, and take the adjoint of the second and multiply by $\alpha$ on the right:
\begin{eqnarray*}
	\beta \circ f \circ \alpha	& = &	\beta \circ g_\alpha \circ (1 - \tau)	\\
	\beta \circ f \circ \alpha	& = &	(1 - \tau^{-1}) \circ h_\beta^\ast \circ \alpha	\\
					& = &	(1-\tau) \circ (-\tau^{-1}) \circ h_\beta^\ast \circ \alpha.
\end{eqnarray*}
Equate and apply Lemma~\ref{oneone} to get
\begin{eqnarray*}
	(\beta \circ g_\alpha)^\ast(1)	& = &	((-\tau^{-1}) \circ h_\beta^\ast \circ \alpha)(1)	\\
	g_\alpha^\ast(\beta)		& = &	- \tau^{-1} (h_\beta^\ast(\alpha))	\\
	\langle \alpha, \beta \rangle_f	& = &	- \tau^{-1}(\langle \beta,\alpha \rangle_{f^\ast})	\\
					& = &	- \langle \beta,\alpha \rangle_{f^\ast},
\end{eqnarray*}
since $\langle \beta,\alpha \rangle_{f^\ast} \in \Fq$.
\end{pf}

\begin{theorem}
\label{duality}
If $f \in \FF$ is nonzero, then the pairing
	$$\langle \;\;,\;\; \rangle_f : \ker f \times \ker f^\ast \rightarrow \Fq$$
exhibits $\ker f$ as the Pontryagin dual of $\ker f^\ast$ as a topological $\Fq$-vector space.
\end{theorem}

\begin{pf}
First we check that the pairing is continuous.
Because of Proposition~\ref{backwards}, it will suffice to show that given any bounded subset $B$ of $\ker f^\ast$, there exists a neighborhood $U$ of 0 in $\ker f$ such that $\langle \alpha,\beta \rangle_f=0$ for $\alpha \in U$, $\beta \in B$.
Given $\alpha \in \ker f$, write
	$$f \circ \alpha = g_\alpha \circ (1 - \tau)$$
as in the definition of the pairing.
As $\alpha$ tends to zero, each coefficient of $f \circ \alpha$ tends to zero, and so does each coefficient of $g_\alpha$, by the construction in the proof of Lemma~\ref{remainder}.
Then the same is true for the coefficients of $g_\alpha^\ast$.
Combined with the knowledge that $g_\alpha^\ast$ is everywhere convergent for any $\alpha$, this implies that $g_\alpha^\ast(\beta)$ tends to zero uniformly for $\beta \in B$ as $\alpha$ tends to zero, as desired.
Thus the pairing is continuous.

The map induced by the pairing is
\begin{eqnarray*}
	\Psi: \ker f	& \rightarrow	&	\Hom_\Fq(\ker f^\ast,\Fq)	\\
		 \alpha	& \mapsto	&	g_\alpha^\ast|_{\ker f^\ast}
\end{eqnarray*}
This is well-defined and continuous, since the pairing is.

We now show $\Psi$ is surjective.
If $\phi \in \Hom_\Fq(\ker f^\ast,\Fq)$ is nonzero, then $H=\ker \phi$ is some open and closed subspace in $G=\ker f^\ast$, so $H \in \GG$ by~(6) in Proposition~\ref{class}.
Thus by Theorem~\ref{existence}, there exists $h \in \FF$ such that $\ker h = H$.
Since $H$ has codimension 1 in $G$ (as an $\Fq$-vector space), $h$ must map $G$ to a one-dimensional $\Fq$-vector space inside $C$.
By scaling $h$, we may assume $h(G)=\Fq$.
Now $(1 - \tau) \circ h$ has kernel $G$, so by the uniqueness in Theorem~\ref{existence},
	$$\tau^n \alpha \circ f^\ast = (1 - \tau) \circ h$$
for some $\alpha \in C^\ast$ and $n \in \Z$.
Take adjoints and multiply on the right by $\tau^n$:
\begin{eqnarray*}
	f \circ \alpha	& = &	h^\ast \circ (1-\tau^{-1}) \circ \tau^n	\\
			& = &	h^\ast \circ (-\tau^{n-1}) \circ (1-\tau).
\end{eqnarray*}
Evaluating at 1 shows that $\alpha \in \ker f$.
Comparing with~(\ref{def}) shows
\begin{eqnarray*}
	g_\alpha	& = &	h^\ast \circ (-\tau^{n-1})	\\
	g_\alpha^\ast	& = &	- \tau^{1-n} \circ h		\\
	g_\alpha^\ast(G) & = &	- \tau^{1-n}(h(G))		\\
			& = &	\tau^{1-n}(\Fq)			\\
			& = & \Fq.
\end{eqnarray*}
Thus $\phi$ and $\Psi(\alpha)=g_\alpha^\ast|_{\ker f^\ast}$ are nonzero $\Fq$-homomorphisms from $G$ to $\Fq$ with the same kernel $H$, so $\phi = c \Psi(\alpha) = \Psi(c \alpha)$ for some $c \in \Fq^\ast$.
Therefore $\Psi$ is surjective.

Applying the above argument to $f^\ast$ and invoking Proposition~\ref{backwards} shows that the other induced map
	$$\ker f^\ast \rightarrow \Hom_\Fq(\ker f,\Fq)$$
is surjective as well.
Now apply Lemma~\ref{showdual3}.
\end{pf}

\begin{rems}
If $f = \sum_{i=m}^n a_i \tau^i$ is a finite series with $a_i$ in {\em any} field $L$ of characteristic $p$, then the construction yields a perfect pairing between the kernels of $f$ and $f^\ast$ acting on $\Lbar$.
Here, in addition, the pairing will be $\Gal(L\sep/L)$-equivariant, since its construction is canonical.
The same will holds for infinite series, when $L$ has a non-archimedean valuation and $L\sep$ is replaced by topological Galois closure.
\end{rems}

\medskip

Next we prove a compatibility result.

\begin{prop}
\label{compatible}
Let $f$ and $h$ be nonzero elements of $\FF$.
Then for all $\alpha \in \ker(f\circ h)$ and $\beta \in \ker f^\ast \subseteq \ker (f \circ h)^\ast$,
	$$\langle \alpha,\beta \rangle_{f \circ h} = \langle h(\alpha),\beta \rangle_f.$$
Similarly, for all $\alpha \in \ker f \subseteq \ker(h \circ f), \beta \in \ker(h \circ f)^\ast$,
	$$\langle \alpha,\beta \rangle_{h \circ f} = \langle \alpha,h^\ast(\beta) \rangle_f.$$
\end{prop}

\begin{pf}
For the first part, write
\begin{eqnarray*}
	(f \circ h) \circ \alpha		& = &	G_\alpha \circ (1 - \tau)	\\
	f \circ h(\alpha)	& = &	g_{h(\alpha)} \circ (1 - \tau)
\end{eqnarray*}
Subtract to get
	$$f \circ (h \circ \alpha - h(\alpha)) = (G_\alpha-g_{h(\alpha)}) \circ (1 - \tau).$$
Since $h \circ \alpha - h(\alpha)$ kills $1$,
	$$h \circ \alpha - h(\alpha) = q \circ (1 - \tau)$$
for some $q \in \FF$, by Lemma~\ref{remainder}.
Substitute and cancel $1-\tau$ on the right:
	$$f \circ q = G_\alpha - g_{h(\alpha)}.$$
Take adjoints and apply both sides to $\beta \in \ker f^\ast$ to get
\begin{eqnarray*}
	0	& = &	G_\alpha^\ast(\beta) - g_{h(\alpha)}^\ast(\beta)	\\
		& = &	\langle \alpha, \beta \rangle_{f \circ h} - \langle h(\alpha),\beta \rangle_f
\end{eqnarray*}
as desired.

For the second part, we have
\begin{eqnarray*}
	\langle \alpha,\beta \rangle_{h \circ f}
		& = &	- \langle \beta, \alpha \rangle_{f^\ast \circ h^\ast} \;\;\; \text{(by Proposition~\ref{backwards})}	\\
		& = &	- \langle h^\ast(\beta), \alpha \rangle_{f^\ast} \;\;\; \text{(by what we just showed)}	\\
		& = &	\langle \alpha,h^\ast(\beta) \rangle_f,
\end{eqnarray*}
by Proposition~\ref{backwards} again.
\end{pf}

\begin{prop}
\label{annihilator}
Suppose $h,g \in \FF$ and $f=h \circ g$.
Then
	$$\{\, \beta \in \ker f^\ast : \langle \alpha,\beta \rangle_f = 0 \;\; \forall \alpha \in \ker g \,\} = \ker h^\ast.$$
\end{prop}

\begin{pf}
Suppose $\beta \in \ker f^\ast$.
Then, by Proposition~\ref{compatible},
\begin{eqnarray*}
	\langle \alpha,\beta \rangle_f = 0 \;\; \forall \alpha \in \ker g
		& \iff & \langle \alpha,h^\ast(\beta) \rangle_g = 0 \;\; \forall \alpha \in \ker g	\\
		& \iff & h^\ast(\beta)=0 \;\;\; \text{ (by Theorem~\ref{duality})}.
\end{eqnarray*}
\end{pf}

Below are three more properties of the pairing.
These were discovered by Elkies~\cite{elkies} (for the case of additive polynomials), although some of our proofs are new.
We will not need these in the sequel, but they are of interest in their own right.
First we explain how the pairing changes if we change the underlying finite field.

\begin{prop}
Let $q'=q^m$, for some $m \ge 1$.
If $f \in \FF$ is actually $\Fqq$-linear, then we can define another pairing $\langle \;\;,\;\; \rangle'_f$ between $\ker f$ and $\ker f^\ast$, this time taking values in $\Fqq$.
(Note that $f^\ast$ is the same whether $f$ is considered as $\Fq$-linear or $\Fqq$-linear.)
Then $\langle \alpha,\beta \rangle_f = \Tr_{\Fqq/\Fq} \langle \alpha,\beta \rangle'_f$.
\end{prop}

\begin{pf}
Write	
\begin{eqnarray*}
	f \circ \alpha	& = &	G_\alpha \circ (1 - \tau^m)	\\
			& = &	G_\alpha \circ (1+ \tau + \cdots + \tau^{m-1}) \circ (1-\tau).
\end{eqnarray*}
By definition of the pairing,
\begin{eqnarray*}
	\langle \alpha,\beta \rangle_f	& = & (G_\alpha \circ (1+ \tau + \cdots + \tau^{m-1}))^\ast(\beta)	\\
		& = &	(1+\tau + \cdots + \tau^{m-1})^\ast(G_\alpha^\ast(\beta))	\\
		& = &	(1+\tau^{-1} + \cdots + \tau^{-(m-1)})(\langle \alpha,\beta \rangle'_f)	\\
		& = &	\Tr_{\Fqq/\Fq} \langle \alpha,\beta \rangle'_f,	\\
\end{eqnarray*}
since $\Gal(\Fqq/\Fq) = \{ 1, \tau^{-1},\ldots,\tau^{-(m-1)} \}$.
\end{pf}

\begin{theorem}
\label{skewsymmetric}
Suppose $f \in \FF$ is nonzero and $f=f^\ast$.
Then
\begin{enumerate}
	\item $\langle \alpha,\alpha \rangle_f = 0$ for all $\alpha \in \ker f$.
	\item If $h \in \FF$ and $f=h^\ast \circ h$, then $\ker h$ is a maximal isotropic closed sub-$\Fq$-vector space of $\ker f$ under $\langle \;,\; \rangle_f$.
	\item Conversely, if $H$ is a maximal isotropic closed sub-$\Fq$-vector space of $\ker f$ under $\langle \;,\; \rangle_f$, then there exists $h \in \FF$ such that $f=h^\ast \circ h$ and $\ker h = H$.
\end{enumerate}
\end{theorem}

\begin{pf}
If $p \not=2$, then part~(1) follows already from Proposition~\ref{backwards}.
Hence assume $p=2$.
Write $f \circ \alpha = g_\alpha \circ (1-\tau)$ as usual.
Multiply by $\alpha$ on the left, and apply Lemma~\ref{remainder2} to write
\begin{eqnarray*}
	\alpha \circ f \circ \alpha	& = &	[(1-\tau) \circ q + (\alpha \circ g_\alpha)^\ast(1)] \circ (1-\tau)	\\
		& = &	(1-\tau^{-1}) \circ j \circ (1-\tau) + g_\alpha^\ast(\alpha) \circ (1-\tau)	\;\;\;	\text{(where $j=\tau \circ q$).}
\end{eqnarray*}
Subtract this from its adjoint and use $f=f^\ast$ and $g_\alpha^\ast(\alpha) \in \Fq$ to get
	$$0 = (1-\tau^{-1}) \circ (j^\ast-j) \circ (1-\tau) + g_\alpha^\ast(\alpha) \circ (\tau-\tau^{-1}).$$
Since $\tau-\tau^{-1}=(1-\tau^{-1}) \circ (1-\tau)$ in characteristic 2, we can cancel $1-\tau^{-1}$ on the left and $1-\tau$ on the right to obtain
	$$0 = (j^\ast-j) + g_\alpha^\ast(\alpha).$$
Take coefficients of $\tau^0$ to deduce
	$$\langle \alpha,\alpha \rangle_f = g_\alpha^\ast(\alpha) = 0.$$
This proves part~(1).

Now, for part~(2), $h$ is continuous, so $\ker h$ is closed.
Proposition~\ref{annihilator} says that
	$$\{\, \beta \in \ker f^\ast : \langle \alpha,\beta \rangle_f = 0 \;\; \forall \alpha \in \ker h \,\} = \ker h,$$
which proves part~(2).

Finally let us prove~(3).
Assume $H$ is a maximal isotropic closed sub-$\Fq$-vector space of $\ker f$ under $\langle \;,\; \rangle$.
Then $H \in \GG$ by~(6) in Proposition~\ref{class}.
By Theorem~\ref{existence} we can find $g \in \FF$ with kernel $H$, and by Theorem~\ref{factors}, we have $f = j \circ g$ for some $j \in \FF$.
By Proposition~\ref{annihilator},
	$$\ker j^\ast = \{\, \beta \in \ker f^\ast : \langle \alpha,\beta \rangle_f = 0 \;\; \forall \alpha \in H \,\} =H,$$
since $H$ is maximal isotropic.
By the uniqueness in Theorem~\ref{existence}, we have $j^\ast = \epsilon \circ \tau^n \circ g$ where $\epsilon \in C^\ast$ and $n \in \Z$.
Then
	$$f = j \circ g = g^\ast \circ \tau^{-n} \circ \epsilon \circ g.$$
Taking adjoints and using $f=f^\ast$ shows that
	$$f = g^\ast \circ \epsilon \circ \tau^n \circ g.$$
If we equate and cancel $g^\ast$ on the left and $g$ on the right, we obtain
	$$\tau^{-n} \circ \epsilon = \epsilon \circ \tau^n,$$
which forces $n=0$.  Now simply take $h = \sqrt{\epsilon} \circ g$ where $\sqrt{\epsilon}$ is any square root of $\epsilon$ in $C$.
\end{pf}

\begin{theorem}
\label{symmetric}
If $f \in \FF$ and $f^\ast=-f$, then $\langle \alpha,\beta \rangle_f = \langle \beta,\alpha \rangle_f$ for all $\alpha,\beta \in \ker f$.
\end{theorem}

\begin{pf}
By definition of the pairing, $\langle \;,\; \rangle_{-f} = -\langle \;,\; \rangle_f$, so
\begin{eqnarray*}
	\langle \alpha,\beta \rangle_f	& = &	-\langle \beta,\alpha \rangle_{f^\ast}	\quad \text{(by Proposition~\ref{backwards})}	\\
		& = & -\langle \beta,\alpha \rangle_{-f}	\\
		& = & \langle \beta,\alpha \rangle_f.
\end{eqnarray*}
\end{pf}

\begin{rems}
(I thank Noam Elkies for these.)
The theory of additive polynomials shares much with the theory of differential operators, as pointed out by Ore~\cite{Or1}.
In fact the pairing $\langle \;,\; \rangle_f$ can be obtained an analogue of a known pairing on kernels of differential operators.
If $f$ is a differential operator and $f^\ast$ denotes its formal adjoint, then we have an identity $uf^\ast(v)-vf(u)=B'$ where $B$ is a bilinear form in $\{u,u',u'',\ldots,u^{(n-1)}\}$ and $\{v,v',v'',\ldots,v^{(n-1)}\}$.
(See p.\ 124 in~\cite{ince}, where $B$ is called a ``bilinear concomitant.'')
If $u \in \ker f$ and $v \in \ker f^\ast$, then $B(u,v)$ is a constant depending bilinearly on $u$ and $v$.

Similarly, if $f \in \FF$, then one can write $uf^\ast(v)-vf(u)=B^q-B$ where $B$ is a (potentially infinite) linear combination of terms of the form $u^{q^m} v^{q^n}$ ($m,n \in \Z$).
One can check that $B$ converges to a continuous bilinear function on $C \times C$ by writing $B$ explicitly in terms of the coefficients of $f$, and using Proposition~\ref{convergence}.
Clearly $B(\alpha,\beta) \in \Fq$ if $\alpha \in \ker f$ and $\beta \in \ker f^\ast$.
\end{rems}

\begin{prop}
\label{ince}
If $\alpha \in \ker f$ and $\beta \in \ker f^\ast$, then $B(\alpha,\beta)=\langle \alpha,\beta \rangle_f$.
\end{prop}

\begin{pf}
Let $b_\alpha(v)=B(\alpha,v)$.
If we set $u=\alpha$ in the defining equation for $B$, we have as functions of $v$,
\begin{eqnarray*}
	\alpha \circ f^\ast	& = &	(\tau-1) \circ b_\alpha, \\
	f \circ \alpha		& = &	b_\alpha^\ast \circ (\tau^{-1}-1) \\
				& = &	b_\alpha^\ast \circ \tau^{-1} \circ (1 -\tau).
\end{eqnarray*}
Comparing with equation~(\ref{def}) shows that $g_\alpha = b_\alpha^\ast \circ \tau^{-1}$, so
\begin{eqnarray*}
	\langle \alpha,\beta \rangle_f	& = &	g_\alpha^\ast(\beta)	\\
					& = &	(\tau \circ b_\alpha)(\beta)	\\
					& = &	\tau(B(\alpha,\beta))	\\
					& = &	B(\alpha,\beta),
\end{eqnarray*}
since $B(\alpha,\beta) \in \Fq$.
\end{pf}

This alternative definition of the pairing could clearly be used to give new proofs of the properties of the pairing.
For instance, Proposition~\ref{backwards} would be immediate, as would the first part of Theorem~\ref{skewsymmetric}.

\section{$A$-module pairings for Drinfeld modules}
\label{amodule}

Let us retain the assumption that $A$ is the affine ring of a nonsingular projective curve over $\Fq$ minus a closed point $\infty$, and suppose we have an $\Fq$-algebra homomorphism $\iota: A \rightarrow C$.
The kernel $\pp$ of $\iota$ will be called the characteristic.

Let $a \mapsto \phi_a$ and $a \mapsto \psi_a$ be ring homomorphisms from $A$ to $\FF$.
For example, $\phi$ and $\psi$ might be Drinfeld $A$-modules over $C$.
We say that a nonzero $f \in \FF$ is an {\em isogeny} from $\phi$ to $\psi$ if $\psi_a \circ f = f \circ \phi_a$ for all $a \in A$.
If $f: \phi \rightarrow \psi$ is an isogeny, then $\ker f$ is an $A$-module via $\phi$ and $\ker f^\ast$ is an $A$-module via $\psi^\ast$.

\begin{prop}
\label{apairing}
Let $a \mapsto \phi_a$ and $a \mapsto \psi_a$ be ring homomorphisms from $A$ to $\FF$.
Suppose $f \in \FF$ is an isogeny from $\phi$ to $\psi$.
Then for all $a \in A$, $\alpha \in \ker f$, $\beta \in \ker f^\ast$, we have
	$$\langle \phi_a(\alpha),\beta \rangle_f = \langle \alpha,\psi_a^\ast(\beta) \rangle_f,$$
and
\begin{eqnarray*}
	[ \;\; , \;\; ]_f : \ker f \times \ker f^\ast	& \rightarrow &	\Hom_\Fq(A,\Fq)	\\
	\alpha \; , \; \beta \;\;\;\;\;\;\;		& \mapsto &	(a \mapsto \langle \phi_a(\alpha),\beta \rangle)
\end{eqnarray*}
is a pairing of $A$-modules which exhibits $\ker f^\ast$ as the Pontryagin dual of $\ker f$ as an $A$-module.

Furthermore, if $\Psi:A \rightarrow \FF$ is a third ring homomorphism, and $h:\Psi \rightarrow \phi$ is another isogeny, then we have the following two compatibility relations:
if $\alpha \in \ker(f\circ h),\beta \in \ker f^\ast \subseteq \ker (f \circ h)^\ast$, then
	$$[ \alpha,\beta ]_{f \circ h} = [ h(\alpha),\beta ]_f.$$

Similarly, for all $\alpha \in \ker h \subseteq \ker(f \circ h), \beta \in \ker(f \circ h)^\ast$,
	$$[ \alpha,\beta ]_{f \circ h} = [ \alpha,f^\ast(\beta) ]_h.$$
\end{prop}

\begin{pf}
Given $a \in A$, $\alpha \in \ker f$, $\beta \in \ker f^\ast$, we have
\begin{eqnarray*}
	\langle \phi_a(\alpha),\beta \rangle_f	& = &	\langle \alpha,\beta \rangle_{f \circ \phi_a}	\;\;\; \text{(by Proposition~\ref{compatible})}	\\
						& = &	\langle \alpha,\beta \rangle_{\psi_a \circ f}	\\
						& = &	\langle \alpha,\psi_a^\ast(\beta) \rangle_f	\;\;\; \text{(by Proposition~\ref{compatible} again}).
\end{eqnarray*}
Apply Lemma~\ref{showdual3} to get the duality property of $[\;\;,\;\;]_f$.
The compatibility relations follow from those for $\langle \;\;,\;\; \rangle$ proved in Proposition~\ref{compatible}.
\end{pf}

\begin{cor}
\label{weil}
Let $\phi$ be a Drinfeld $A$-module over a field $L$.
Fix a nonzero $a \in A$.
Let $\phi[a]$ be the kernel of $\phi_a$ on $\Lbar$, and similarly define $\phi^\ast[a]=\ker \phi_a^\ast$.
There is a Galois-equivariant perfect pairing of finite $A$-modules
	$$[ \;\;,\;\; ]_a: \phi[a] \times \phi^\ast[a] \rightarrow \widehat{(A/a)} \subset \hat{A}.$$
(Here $\Gal(L\sep/L)$ acts trivially on $\hat{A}$ and $\widehat{(A/a)}$.)
If we also have a nonzero $b \in A$, then for any $\alpha \in \phi[ab]$, $\beta \in \psi^\ast[a]$
	$$[ \alpha,\beta ]_{ab} = [\phi_b(\alpha),\beta]_a.$$
Similarly, for $\alpha \in \phi[a]$, $\beta \in \psi^\ast[ab]$,
	$$[ \alpha,\beta ]_{ab} = [\alpha,\psi^\ast_b(\beta)]_a.$$
\end{cor}

The properties above of $[\;\;,\;\;]_a$ should remind one of the Weil pairing for abelian varieties.
Another proof of the Galois-equivariant duality between $\phi[a]$ and $\phi^\ast[a]$ was discovered independently by Taguchi.
(See the appendix of~\cite{gossadjoint}.)

As another application of our pairings, we can describe the kernel of the adjoint of the exponential function associated with a Drinfeld module.

\begin{cor}
\label{estar}
Let $C$ be the completion of the algebraic closure of $K_\infty$ with $|\;|_\infty$, and let $\phi$ be a Drinfeld $A$-module over $C$.
Let $e(z) \in C[[z]]$ be the associated exponential function.
Then there is a natural pairing
	$$\ker e \times \ker e^\ast \rightarrow \hat{A}$$
which exhibits $\ker e^\ast$ as the Pontryagin dual of $\ker e$ as an $A$-module.
If the lattice $\ker e$ is isomorphic to the direct sum of fractional ideals $I_1 \oplus \cdots \oplus I_r$, then $\ker e^\ast$ is isomorphic to $K_\infty/(I_1^{-1}J) \oplus \cdots \oplus K_\infty/(I_r^{-1}J)$ as a topological $A$-module.
(Here $J$ is as in Theorem~\ref{fractionaldual}.)
\end{cor}

\begin{pf}
Simply note that $e$ is an isogeny from $C$ with the standard $A$-module structure to $\phi$, and apply Proposition~\ref{apairing}.
The last assertion follows from Theorem~\ref{fractionaldual}.
\end{pf}

\begin{prop}
Let $\phi$ and $e$ be as in the previous corollary.
Then $\ker e^\ast$ is the closure of $\bigcup_{a \in A} \phi^\ast[a]$.
\end{prop}

\begin{pf}
Taking adjoints of
	$$e \circ a = \phi_a \circ e$$
yields
	$$a \circ e^\ast = e^\ast \circ \phi^\ast_a,$$
from which it is clear that $\phi^\ast[a] \subset \ker e^\ast$.

Corollary~\ref{estar} says that $\ker e^\ast$ as a topological $A$-module via $\phi^\ast$ is isomorphic to $K_\infty/C_1 \oplus \cdots \oplus K_\infty/C_r$, for some fractional $A$-ideals $C_1,\ldots,C_r$.
The torsion submodule in this module is $K/C_1 \oplus \cdots \oplus K/C_r$, which is dense, so the result follows.
\end{pf}

\section{Tate module pairings}

Because of the compatibility relations, we can use our pairings to construct a pairing between Tate modules, just as the Weil pairing gives rise to a pairing between Tate modules.
Let $\qq$ be a nonzero prime of $A$ different from the characteristic $\pp$ of the Drinfeld module.
Let $K$ be the fraction field of $A$, and let $A_\qq$, $K_\qq$ be the completions at the prime $\qq$.
The Tate module of a Drinfeld $A$-module over $L$ is
	$$T_\qq(\phi) = \Hom_A(K_\qq/A_\qq,\phi(\Lbar))$$
where $\phi(\Lbar)$ denotes the additive group of $\Lbar$ with the $A$-module structure given by $\phi$.
Similarly we define the Tate module of $\phi^\ast$ as
	$$T_\qq(\phi^\ast) = \Hom_A(K_\qq/A_\qq,\phi^\ast(\Lbar)).$$
These are both free $A_\qq$-modules of rank equal to the rank of the Drinfeld module.

\begin{prop}
\label{tate}
There is a Galois-equivariant continuous perfect pairing of $A_\qq$-modules
	$$[\;\;,\;\;]_\qq : T_\qq(\phi) \times T_\qq(\phi^\ast) \rightarrow \Hom_\Fq(K_\qq/A_\qq,\Fq) \cong A_\qq$$.
\end{prop}

\begin{pf}
Given $\alpha \in T_\qq(\phi)$, $\beta \in T_\qq(\phi^\ast)$, define $[\alpha,\beta]_\qq \in \Hom_\Fq(K_\qq/A_\qq,\Fq)$ by
	$$[\alpha,\beta]_\qq(b)=[\alpha(b),\beta(a^{-1})]_a$$
where $a \in A$ kills $b$ and generates an ideal power of $\qq$.
As is well known, $A$ is a Dedekind domain with finite class number, so some power $\qq^h$ of $\qq$ is principal, so it is always possible to find such an $a$, given $b$.
The definition is independent of the choice of $a$, by the last compatibility relation in Corollary~\ref{weil}.

If we fix $a \in A$ generating a power of $\qq$, we can also obtain the pairing as the inverse limit of the pairings
	$$[\;\;,\;\;]_{a^n} : \phi[a^n] \times \phi^\ast[a^n] \longrightarrow \Hom_\Fq(A/a^n,\Fq) = \Hom_\Fq(a^{-n}A/A,\Fq),$$
with respect to the maps
\begin{eqnarray*}
	\phi[a^{n+1}]		& \stackrel{a}{\longrightarrow}	& \phi[a^n]	\\
	\phi^\ast[a^{n+1}]	& \stackrel{a}{\longrightarrow}	& \phi^\ast[a^n]	\\
	\Hom_\Fq(a^{-(n+1)}A/A,\Fq)	& \stackrel{\text{res}}{\longrightarrow}	& \Hom_\Fq(a^{-n}A/A,\Fq).
\end{eqnarray*}
The Galois-equivariance, continuity, and perfectness of the pairing then follow from the properties of the pairings $[\;\;,\;\;]_{a^n}$.

Finally, $\Hom_\Fq(K_\qq/A_\qq,\Fq)$ is (non-canonically) isomorphic to $A_\qq$ as an $A_\qq$-module, by Theorem~\ref{aqdual}.
\end{pf}

As pointed out by Goss, this duality of Tate modules has the following corollary.

\begin{cor}
\label{semisimple}
Suppose the Drinfeld module $\phi$ is defined over a finite extension $L$ of $K$.
Then $$T_\qq(\phi^\ast) \otimes_{A_\qq} K_\qq$$ is a semi-simple $K_\qq[\Gal(L\sep/L)]$-module.
\end{cor}

\begin{pf}
Combine the previous theorem with the main theorem in~\cite{taguchi}.
\end{pf}

\section{Questions}
\label{questions}

Our results show that in two cohomological realizations (the period lattice and Tate module), the cohomology of the adjoint of a Drinfeld module is dual to the cohomology of the original Drinfeld module.
There is a third cohomological realization of a Drinfeld module, namely the de Rham cohomology~\cite{gekeler} developed by Anderson, Deligne, Gekeler, and Yu.

\begin{question}
Is it possible to give a reasonable definition for the de Rham cohomology of the adjoint of a Drinfeld module and prove that it is dual in some sense to the de Rham cohomology of the original Drinfeld module?
\end{question}

It would also be nice to generalize the applications of this paper to include $t$-modules, the higher dimensional analogues of Drinfeld modules.

\begin{question}
Can one prove results similar to those in this paper for fractional power series in more than one variable?
\end{question}

A reasonable approach is to work with the ring $\Md(\FF)$ of $d$-by-$d$ matrices with coefficients in $\FF$.
Elements of this ring act on $C^d$ in an obvious way.
Moreover, there is a multiplication-reversing adjoint map $\Md(\FF) \rightarrow \Md(\FF)$, which maps a matrix $A=(a_{ij})$ to $A^\ast=(a_{ji}^\ast)$, i.e., the transpose of the matrix obtained by taking the adjoint of each entry.
Hence, for example, it is possible to define the adjoint of a $t$-module
\begin{eqnarray*}
  \phi:	\Fq[t]	& \rightarrow	& \Md(\Ct) \subset \Md(\FF)	\\
	a	& \mapsto	& \phi_a
\end{eqnarray*}
as
\begin{eqnarray*}
  \phi^\ast:	\Fq[t]	& \rightarrow	& \Md(\FF)	\\
		a	& \mapsto	& \phi_a^\ast.
\end{eqnarray*}

For each $A = (a_{ij}) \in \Md(\FF)$ it is also possible to define a natural bilinear pairing
	$$\langle \;,\; \rangle_A : \ker A \times \ker A^\ast \rightarrow \Fq$$
by following the differential operator analogy.
Let $\cdot$ denote the standard inner product on $C^d$.
As in the remarks at the end of Section~\ref{adjointpd}, for each $f \in \FF$, let $B_f(u,v)$ be the bilinear function of $u,v \in C$ such that
	$$u f^\ast(v) - v f(u) = B_f^q - B_f.$$
Then, for $u,v \in C^d$, the bilinear function
	$$B_A(u,v) = \sum_{i,j} B_{a_{ij}}(u_i,v_j)$$
satisfies
	$$u \cdot A^\ast(v) - v \cdot A(u) = B_A^q-B_A,$$
and we may define
	$$\langle u,v \rangle_A = B_A(u,v).$$
Thus for instance, one has a pairing between the $a$-torsion of a $t$-module and the $a$-torsion of its adjoint.

\section*{Acknowledgements}

Many thanks go to David Goss, who introduced me to adjoints of polynomials and Drinfeld modules, and asked some of the questions answered by this paper.
I thank T.\ Y.\ Lam, Hendrik Lenstra, and Ken Ribet for helping me find references for Newton polygons.
Thanks also to Noam Elkies for sharing his preprint with me, and to Michael Rosen for suggesting references of use in proving Theorem~\ref{fractionaldual}.
Finally, I thank the referee for many intelligent and helpful suggestions.


\end{document}